\documentclass[11pt]{article}
\usepackage{amsmath}
\usepackage{amssymb}

\usepackage{theorem}
\numberwithin{equation}{section}

\newtheorem{theorem}{Theorem}[section]
\newtheorem{proposition}[theorem]{Proposition}
\newtheorem{corollary}[theorem]{Corollary}
\newtheorem{conjecture}[theorem]{Conjecture}
\newtheorem{lemma}[theorem]{Lemma}
\newtheorem{definition}{Definition}[section]

\begin{document}
\title{Global well-posedness and scattering for the defocusing, $L^{2}$-critical, nonlinear Schr{\"o}dinger equation when $d = 1$}
\date{\today}
\author{Benjamin Dodson}
\maketitle

\noindent \textbf{Abstract:} In this paper we prove that the defocusing, quintic nonlinear Schr{\"o}dinger initial value problem is globally well-posed and scattering for $u_{0} \in L^{2}(\mathbf{R})$. To do this, we will prove a frequency localized interaction Morawetz estimate similar to the estimate made in \cite{CKSTT4}. Since we are considering an $L^{2}$ - critical initial value problem we will localize to low frequencies.

\section{Introduction} The quintic nonlinear Schr{\"o}dinger initial value problem is given by

\begin{equation}\label{0.1}
\aligned
i u_{t} + \Delta u &= F(u), \\
u(0,x) &= u_{0} \in L^{2}(\mathbf{R}),
\endaligned
\end{equation}

\noindent where $F(u) = \mu |u|^{4} u$, $\mu = \pm 1$, $u(t) : \mathbf{R} \rightarrow \mathbf{C}$. When $\mu = +1$ $(\ref{0.1})$ is said to be defocusing and when $\mu = -1$ $(\ref{0.1})$ is said to be focusing. It was observed in \cite{CaWe} that the solution to $(\ref{0.1})$ conserves mass,

\begin{equation}\label{0.2}
M(u(t)) = \int |u(t,x)|^{2} dx = M(u(0)),
\end{equation}

\noindent and energy

\begin{equation}\label{0.2.1}
E(u(t)) = \frac{1}{2} \int |\nabla u(t,x)|^{2} dx + \frac{\mu}{6} \int |u(t,x)|^{6} dx = E(u(0)).
\end{equation}

\noindent The initial value problem $(\ref{0.1})$ also obeys a scaling symmetry. If $u(t,x)$ is a solution to $(\ref{0.1})$ on a time interval $[0, T]$, then

\begin{equation}\label{0.2.2}
u_{\lambda}(t,x) = \frac{1}{\lambda^{1/2}} u(\frac{t}{\lambda^{2}}, \frac{x}{\lambda})
\end{equation}

\noindent is a solution to $(\ref{0.1})$ on $[0, \lambda^{2} T]$ with $u(0,x) = \frac{1}{\lambda^{1/2}} u_{0}(\frac{x}{\lambda})$.

\begin{equation}\label{0.2.2.1}
 \| \frac{1}{\lambda^{1/2}} u_{0}(\frac{x}{\lambda}) \|_{L^{2}(\mathbf{R})} = \| u_{0}(x) \|_{L^{2}(\mathbf{R})}.
\end{equation}

\noindent Therefore, $(\ref{0.1})$ is called $L^{2}$ - critical or mass critical.\vspace{5mm}

\noindent A solution to $(\ref{0.1})$ obeys Duhamel's formula

\begin{definition}\label{d0.1}
 $u : I \times \mathbf{R}^{d} \rightarrow \mathbf{C}$, $I \subset \mathbf{R}$ is a solution to $(\ref{0.1})$ if for any compact $J \subset I$, $u \in C_{t}^{0} L_{x}^{2}(J \times \mathbf{R}^{d}) \cap L_{t,x}^{\frac{2(d + 2)}{d}}(J \times \mathbf{R}^{d})$, and for all $t, t_{0} \in I$,

\begin{equation}\label{0.3.1}
 u(t) = e^{i(t - t_{0}) \Delta} u(t_{0}) - i \int_{t_{0}}^{t} e^{i(t - \tau) \Delta} F(u(\tau)) d\tau.
\end{equation}

\end{definition}

\noindent The space $L_{t,x}^{6}(J \times \mathbf{R})$ arises from the Strichartz estimates. This norm is also invariant under the scaling $(\ref{0.2.2})$.

\begin{definition}\label{d0.2}
 A solution to $(\ref{0.1})$ defined on $I \subset \mathbf{R}$ blows up forward in time if there exists $t_{0} \in I$ such that 

\begin{equation}\label{0.3.2}
 \int_{t_{0}}^{\sup(I)} \int |u(t,x)|^{6} dx dt = \infty.
\end{equation}

\noindent $u$ blows up backward in time if there exists $t_{0} \in I$ such that

\begin{equation}\label{0.3.3}
 \int_{\inf(I)}^{t_{0}} \int |u(t,x)|^{6} dx dt = \infty.
\end{equation}

\end{definition}

\begin{definition}\label{d0.0.2}
A solution $u(t,x)$ to $(\ref{0.1})$ is said to scatter forward in time if there exists $u_{+} \in L^{2}(\mathbf{R}^{d})$ such that

\begin{equation}\label{0.10}
\lim_{t \rightarrow \infty} \| e^{it \Delta} u_{+} - u(t,x) \|_{L^{2}(\mathbf{R}^{d})} = 0.
\end{equation}

\noindent A solution is said to scatter backward in time if there exists $u_{-} \in L^{2}(\mathbf{R}^{d})$ such that

\begin{equation}\label{0.11}
\lim_{t \rightarrow -\infty} \| e^{it \Delta} u_{-} - u(t,x) \|_{L^{2}(\mathbf{R}^{d})} = 0.
\end{equation}
\end{definition}

\begin{theorem}\label{t0.0.1}
If $\| u_{0} \|_{L^{2}(\mathbf{R})}$ is sufficiently small, then $(\ref{0.1})$ is globally well-posed and scatters to a free solution as $t \rightarrow \pm \infty$.
\end{theorem}

\noindent \emph{Proof:} See \cite{CaWe}, \cite{CaWe1}. $\Box$\vspace{5mm}

\noindent We will recall the proof of this theorem in $\S 2$. \cite{CaWe}, \cite{CaWe1} also proved that $(\ref{0.1})$ is locally well-posed for $u_{0} \in L^{2}(\mathbf{R})$ on some interval $[0, T]$, where $T(u_{0}) > 0$ depends on the profile of the initial data, not just it size $\| u_{0} \|_{L^{2}(\mathbf{R})}$.

\begin{theorem}\label{t0.0.0.1}
 Given $u_{0} \in L^{2}(\mathbf{R}^{2})$ and $t_{0} \in \mathbf{R}$, there exists a maximal lifespan solution $u$ to $(\ref{0.1})$ defined on $I \subset \mathbf{R}$ with $u(t_{0}) = u_{0}$. Moreover,\vspace{5mm}

1. $I$ is an open neighborhood of $t_{0}$.

2. If $\sup(I)$ or $\inf(I)$ is finite, then $u$ blows up in the corresponding time direction.

3. The map that takes initial data to the corresponding solution is uniformly continuous on compact time intervals for bounded sets of initial data.

4. If $\sup(I) = \infty$ and $u$ does not blow up forward in time, then $u$ scatters forward to a free solution. If $\inf(I) = -\infty$ and $u$ does not blow up backward in time, then $u$ scatters backward to a free solution.
\end{theorem}

\noindent \emph{Proof:} See \cite{CaWe}, \cite{CaWe1}. $\Box$\vspace{5mm}

\noindent There are known counterexamples to $(\ref{0.1})$ globally well-posed and scattering in the focusing case, $\mu = -1$. There are no known counterexamples in the defocusing case. Therefore, it has been conjectured

\begin{conjecture}\label{c0.0.2}
For $d \geq 1$, the defocusing, mass critical nonlinear Schr{\"o}dinger initial value problem $(\ref{0.1})$ is globally well-posed for $u_{0} \in L^{2}(\mathbf{R}^{d})$ and all solutions scatter to a free solution as $t \rightarrow \pm \infty$.
\end{conjecture}

\noindent This conjecture has already been verified for $d \geq 2$.

\begin{theorem}\label{t0.1}
When $d = 2$, $(\ref{0.1})$ is globally well-posed and scattering for $u_{0} \in L^{2}(\mathbf{R}^{2})$.
\end{theorem}

\noindent \emph{Proof:} See \cite{KTV} for a proof in the radial case, \cite{D3} for a proof in the non-radial case.

\begin{theorem}\label{t0.1.1}
When $d \geq 3$, $(\ref{0.1})$ is globally well-posed and scattering for $u_{0} \in L^{2}(\mathbf{R}^{d})$.
\end{theorem}

\noindent \emph{Proof:} See \cite{KVZ}, \cite{TVZ2} for a proof in the radial case, \cite{D2} for a proof in the nonradial case.\vspace{5mm}

\noindent In this paper we tackle the case $d = 1$ and prove

\begin{theorem}\label{t0.2}
$(\ref{0.1})$ is globally well-posed and scattering for $u_{0} \in L^{2}(\mathbf{R})$, $\mu = +1$.
\end{theorem}

\noindent This completes the proof of the conjecture in the defocusing case.\vspace{5mm}

\noindent \textbf{Remark:} \cite{KTV} and \cite{KVZ} also proved global well-posedness and scattering for the focusing, mass-critical initial value problem

\begin{equation}\label{0.2.4.1}
\aligned
iu_{t} + \Delta u &= -|u|^{4/d} u,\\
u(0,x) &= u_{0},
\endaligned
\end{equation}

\noindent with radial data and mass less than the mass of the ground state when $d \geq 2$. Much of the analysis in this paper carries over directly to the focusing case. Therefore, whenever possible we will prove theorems without regard for the sign of $\mu$.\vspace{5mm}

\noindent \textbf{Outline of the Proof.} In this paper we use the concentration compactness method, which is a modification of the induction on energy method. The induction on energy method was introduced in \cite{B2} to prove global well-posedness and scattering for the defocusing energy-critical initial value problem in $\mathbf{R}^{3}$ for radial data.\vspace{5mm}

\noindent \cite{KTV}, \cite{KVZ}, \cite{TVZ2}, \cite{D2}, and \cite{D3} used the concentration compactness method. Since $(\ref{0.1})$ is globally well-posed for small $\| u_{0} \|_{L^{2}(\mathbf{R})}$, if $(\ref{0.1})$, $\mu = +1$ is not globally well-posed for all $u_{0} \in L^{2}(\mathbf{R})$, then there must be a minimum $\| u_{0} \|_{L^{2}(\mathbf{R})} = m_{0}$ where global well-posedness fails. \cite{TVZ1} showed that for conjecture $\ref{c0.0.2}$ to fail, there must exist a minimal mass blowup solution with a number of additional properties.

\begin{theorem}\label{t0.5}
Suppose conjecture $\ref{c0.0.2}$ fails when $d = 1$. Then there exists a maximal lifespan solution on $I \subset \mathbf{R}$, $[0, \infty) \subset I$, $\| u(t) \|_{L_{x}^{2}(\mathbf{R}^{d})} = m_{0}$ which is almost periodic modulo scaling and blows up both forward and backward in time. Moreover, $N(t) \leq 1$ for $t \in [0, \infty)$, $N(0) = 1$, and

\begin{equation}\label{2.5.1}
\int_{0}^{\infty} \int |u(t,x)|^{6} dx = \infty.
\end{equation}

\noindent Additionally, there exists a set $K \subset L^{2}(\mathbf{R})$, $K$ is precompact in $L^{2}(\mathbf{R}^{d})$ such that for all $t \in I$ there exists $Q_{t} \in K$, $x(t), \xi(t) : I \rightarrow \mathbf{R}$ with

\begin{equation}\label{2.5.2}
u(t,x) = \frac{1}{N(t)^{1/2}} e^{ix \cdot \xi(t)} Q_{t}(\frac{x - x(t)}{N(t)}).
\end{equation}

\end{theorem}

\noindent \emph{Proof:} See \cite{KTV}, \cite{TVZ1}, and section four of \cite{TVZ2}.\vspace{5mm}

\noindent \textbf{Remark:} This is also true for a minimal mass blowup solution to the focusing problem $(\ref{0.1})$, $\mu = -1$.\vspace{5mm}

\noindent We will then consider two subcases separately,

\begin{equation}\label{2.5.3}
\int_{0}^{\infty} N(t)^{3} dt < \infty,
\end{equation}

\noindent and

\begin{equation}\label{2.5.4}
\int_{0}^{\infty} N(t)^{3} dt = \infty.
\end{equation}

\noindent We will exclude $(\ref{2.5.3})$ by proving additional regularity, which prevents $N(t) \searrow 0$ as $t \rightarrow \infty$. For $(\ref{2.5.4})$ we will not prove any additional regularity. Instead, we will rely on a frequency localized interaction Morawetz estimate. (See \cite{CKSTT4} for such an estimate in the energy-critical case.) Since we are truncating to low frequencies, our method is very similar to the almost Morawetz estimates that are often used in conjunction with the I-method. (See \cite{B1}, \cite{CKSTT1}, \cite{CKSTT2}, \cite{CKSTT3}, \cite{CR}, \cite{CGT}, \cite{D}, \cite{D1}, \cite{DPST}, and \cite{DPST1} for more information on the I-method.)\vspace{5mm}

\section{Function Spaces and linear estimates}
\noindent \textbf{Linear Strichartz Estimates:}

\begin{definition}\label{d2.0}
A pair $(p,q)$ will be called an admissible pair for $d = 1$ if $\frac{2}{p} = (\frac{1}{2} - \frac{1}{q})$, and $p \geq 4$.
\end{definition}

\begin{theorem}\label{t2.0.1}
If $u(t,x)$ solves the initial value problem

\begin{equation}\label{2.0.1}
\aligned
i u_{t} + \Delta u &= F(t), \\
u(0,x) &= u_{0},
\endaligned
\end{equation}

\noindent on an interval $I$, then

\begin{equation}\label{2.0.2}
\| u \|_{L_{t}^{p} L_{x}^{q}(I \times \mathbf{R})} \lesssim_{p,q,\tilde{p},\tilde{q}} \| u_{0} \|_{L^{2}(\mathbf{R})} + \| F \|_{L_{t}^{\tilde{p}'} L_{x}^{\tilde{q}'}(I \times \mathbf{R})},
\end{equation}

\noindent for all admissible pairs $(p,q)$, $(\tilde{p}, \tilde{q})$. $\tilde{p}'$ denotes the Lebesgue dual of $\tilde{p}$.
\end{theorem}

\noindent \emph{Proof:} See \cite{Tao}.\vspace{5mm}

\noindent $(\ref{2.0.2})$ motivates the definition of the Strichartz space.

\begin{definition}\label{d2.0.2}
Define the norm

\begin{equation}\label{2.0.3}
\| u \|_{S^{0}(I \times \mathbf{R})} \equiv \sup_{(p,q) \text{ admissible }} \| u \|_{L_{t}^{p} L_{x}^{q}(I \times \mathbf{R})}.
\end{equation}

\begin{equation}\label{2.0.3.1}
S^{0}(I \times \mathbf{R}) = \{ u : \| u \|_{S^{0}(I \times \mathbf{R})} < \infty \}.
\end{equation}

\noindent We also define the space $N^{0}(I \times \mathbf{R})$ to be the space dual to $S^{0}(I \times \mathbf{R})$ with appropriate norm. Then in fact,

\begin{equation}\label{2.0.4}
\| u \|_{S^{0}(I \times \mathbf{R})} \lesssim \| u_{0} \|_{L^{2}(\mathbf{R})} + \| F \|_{N^{0}(I \times \mathbf{R})}.
\end{equation}
\end{definition}

\begin{theorem}\label{t2.0.1.1}
$(\ref{0.1})$ is globally well-posed when $\| u_{0} \|_{L^{2}(\mathbf{R})}$ is small.
\end{theorem}

\noindent \emph{Proof:} By $(\ref{2.0.4})$ and the definition of $S^{0}$, $N^{0}$,

\begin{equation}\label{2.0.4.1}
\aligned
\| u \|_{S^{0}((-\infty, \infty) \times \mathbf{R})} &\lesssim \| u_{0} \|_{L^{2}(\mathbf{R})} + \| u \|_{L_{t,x}^{6}((-\infty, \infty) \times \mathbf{R})}^{5} \\ &\lesssim \| u_{0} \|_{L^{2}(\mathbf{R})} + \| u \|_{S^{0}((-\infty, \infty) \times \mathbf{R})}^{5}.
\endaligned
\end{equation}

\noindent By the continuity method, if $\| u_{0} \|_{L^{2}(\mathbf{R})}$ is sufficiently small, then we have global well-posedness. We can also obtain scattering with this argument. $\Box$\vspace{5mm}

\noindent Now define the function

\begin{equation}\label{2.0.4.2}
A(m) = \sup \{ \| u \|_{S^{0}((-\infty, \infty) \times \mathbf{R}^{2})} : \text{ u solves $(\ref{0.1})$}, \| u(0) \|_{L^{2}(\mathbf{R}^{2})} = m \}.
\end{equation}

\noindent If we can prove $A(m) < \infty$ for any $m$, then we have proved global well-posedness and scattering.\vspace{5mm}

\noindent Using a stability lemma from \cite{TVZ1} we can prove that $A(m)$ is an upper semicontinuous function of $m$, which proves that $\{ m : A(m) = \infty \}$ is a closed set. This implies that if global well-posedness and scattering does not hold in the defocusing case for all $u_{0} \in L^{2}(\mathbf{R})$, then there must be a minimum $m_{0}$ with $A(m_{0}) = \infty$. We will discuss the properties of a minimal mass blowup solution more in the next section.\vspace{5mm}

\noindent We will also need the Littlewood-Paley decomposition at various points throughout the paper. Let $\phi \in C_{0}^{\infty}(\mathbf{R})$, radial, $0 \leq \phi \leq 1$,

\begin{equation}\label{2.0.4.3}
\phi(x) = \left\{
            \begin{array}{ll}
              1, & \hbox{$|x| \leq 1$;} \\
              0, & \hbox{$|x| > 2$.}
            \end{array}
          \right.
\end{equation}

\noindent Then define the frequency truncation

\begin{equation}\label{2.0.4.4}
\mathcal F(P_{\leq N} u) = \phi(\frac{\xi}{N}) \hat{u}(\xi).
\end{equation}

\noindent Let $P_{> N} u = u - P_{\leq N} u$ and $P_{N} u = P_{\leq 2N} u - P_{\leq N} u$. We will also depart from the customary notation and say

\begin{equation}\label{2.0.4.5}
P_{1/2} u = P_{\leq 1} u.
\end{equation}

\noindent Throughout the paper it will be necessary to make a Littlewood-Paley decomposition with $\xi_{0} \neq 0$ at the origin. Let

\begin{equation}\label{2.0.4.6}
\tilde{P}_{N, \xi_{0}} u = e^{ix \cdot \xi_{0}} P_{N} (e^{-ix \cdot \xi_{0}} u).
\end{equation}

\noindent \textbf{Function Spaces}

\noindent We utilize the function spaces which are a superposition of free solutions to the Schrodinger equation. See \cite{KoTa}, \cite{HHK} for more information.\vspace{5mm}

\begin{definition}\label{d2.1}
Let $1 \leq p < \infty$. Then $U_{\Delta}^{p}$ is an atomic space, where atoms are piecewise solutions to the linear equation.

\begin{equation}\label{2.1}
u = \sum_{k} 1_{[t_{k}, t_{k + 1})} e^{it \Delta} u_{k}, \hspace{5mm} \sum_{k} \| u_{k} \|_{L^{2}}^{p} = 1.
\end{equation}

\noindent For any function $u$,

\begin{equation}\label{2.1.1}
\| u \|_{U_{\Delta}^{p}} = \inf \{ \sum_{\lambda} |c_{\lambda}| : u = \sum_{\lambda} c_{\lambda} u_{\lambda}, \text{$u_{\lambda}$ are $U_{\Delta}^{p}$ atoms} \}
\end{equation}

\end{definition}

\noindent For any $1 \leq p < \infty$, $U_{\Delta}^{p} \subset L^{\infty} L^{2}$. Additionally, $U_{\Delta}^{p}$ functions are continuous except at countably many points and right continuous everywhere.

\begin{definition}\label{d2.2}
Let $1 \leq p < \infty$. Then $V_{\Delta}^{p}$ is the space of right continuous functions $u \in L^{\infty}(L^{2})$ such that

\begin{equation}\label{2.2}
\| v \|_{V_{\Delta}^{p}}^{p} = \| v \|_{L^{\infty}(L^{2})}^{p} + \sup_{\{ t_{k} \} \nearrow} \sum_{k} \| e^{-it_{k} \Delta} v(t_{k}) - e^{-it_{k + 1} \Delta} v(t_{k + 1}) \|_{L^{2}}^{p}.
\end{equation}

\noindent The supremum is taken over increasing sequences $t_{k}$.
\end{definition}

\begin{theorem}\label{t2.3}
The function spaces $U_{\Delta}^{p}$, $V_{\Delta}^{q}$ obey the embeddings

\begin{equation}\label{2.3}
U_{\Delta}^{p} \subset V_{\Delta}^{p} \subset U_{\Delta}^{q} \subset L^{\infty} (L^{2}), \hspace{5mm} p < q.
\end{equation}

\noindent Let $DU_{\Delta}^{p}$ be the space of functions

\begin{equation}\label{2.4}
DU_{\Delta}^{p} = \{ (i \partial_{t} + \Delta)u ; u \in U_{\Delta}^{p} \}.
\end{equation}

\noindent There is the easy estimate

\begin{equation}\label{2.5}
\| u \|_{U_{\Delta}^{p}} \lesssim \| u(0) \|_{L^{2}} + \| (i \partial_{t} + \partial_{x}^{2}) u \|_{DU_{\Delta}^{p}}.
\end{equation}

\noindent Finally, there is the duality relation

\begin{equation}\label{2.6}
(DU_{\Delta}^{p})^{\ast} = V_{\Delta}^{p'}.
\end{equation}

\noindent These spaces are also closed under truncation in time.

\begin{equation}\label{2.7}
\aligned
\chi_{I} : U_{\Delta}^{p} \rightarrow U_{\Delta}^{p}, \\
\chi_{I} : V_{\Delta}^{p} \rightarrow V_{\Delta}^{p}.
\endaligned
\end{equation}
\end{theorem}

\noindent \emph{Proof:} See \cite{HHK}. $\Box$

\begin{lemma}\label{l2.1}
Suppose $J = I_{1} \cup I_{2}$, $I_{1} = [a, b]$, $I_{2} = [b, c]$, $a \leq b \leq c$.

\begin{equation}\label{2.7.1}
\aligned
\| u \|_{U_{\Delta}^{p}(J \times \mathbf{R})}^{p} \leq \| u \|_{U_{\Delta}^{p}(I_{1} \times \mathbf{R})}^{p} + \| u \|_{U_{\Delta}^{p}(I_{2} \times \mathbf{R})}^{p} \\
\| u \|_{U_{\Delta}^{p}(I_{1} \times \mathbf{R})} \leq \| u \|_{U_{\Delta}^{p}(J \times \mathbf{R}^{d})}.
\endaligned
\end{equation}
\end{lemma}

\noindent \emph{Proof:} See \cite{D3}.\vspace{5mm}

\begin{proposition}\label{p2.2}
\begin{equation}\label{2.7.3}
\| P_{N}((e^{it \Delta} u_{0})(e^{-it \Delta} v_{0})) \|_{L_{t,x}^{2}(\mathbf{R} \times \mathbf{R})} \lesssim \frac{1}{N^{1/2}} \| u_{0} \|_{L^{2}(\mathbf{R})} \| v_{0} \|_{L^{2}(\mathbf{R})}.
\end{equation}

\noindent If the supports of $\hat{u}_{0}(\xi)$ and $\hat{v}_{0}(\xi)$ are separated by distance $N$,

\begin{equation}\label{2.7.2}
\| (e^{it \Delta} u_{0})(e^{it \Delta} v_{0}) \|_{L_{t,x}^{2}(\mathbf{R} \times \mathbf{R})} \lesssim \frac{1}{N^{1/2}} \| u_{0} \|_{L^{2}(\mathbf{R})} \| v_{0} \|_{L^{2}(\mathbf{R})}.
\end{equation}
\end{proposition}

\noindent \emph{Proof:} We prove $(\ref{2.7.3})$. $$\tilde{G}(\tau, \xi) = \int e^{-it \tau} \int_{\xi = \eta_{1} + \eta_{2}} e^{-it \eta_{1}^{2}} e^{it \eta_{2}^{2}} \hat{u}_{0}(\eta_{1}) \hat{v}_{0}(\eta_{2}) d\eta_{1} dt$$

$$ = \int_{\xi = \eta_{1} + \eta_{2}} \delta(\tau + \eta_{1}^{2} - \eta_{2}^{2}) \hat{u}_{0}(\eta_{1}) \hat{v}_{0}(\eta_{2}) d\eta_{1}.$$

\noindent Take $\tilde{F}(\tau, \xi)$ with $\| \tilde{F}(\tau, \xi) \|_{L_{\tau, \xi}^{2}(\mathbf{R} \times \mathbf{R})} = 1$, $F$ supported on $|\xi| \sim N$.

$$\int \int \tilde{F}(-\tau, -\xi) \tilde{G}(\tau, \xi) d\tau d\xi = \int \int \tilde{F}((\eta_{1} + \eta_{2})(\eta_{1} - \eta_{2}), \eta_{1} + \eta_{2}) \hat{u}_{0}(\eta_{1}) \hat{v}_{0}(\eta_{2}) d\eta_{1} d\eta_{2}.$$

\noindent Making a change of variables, this proves $(\ref{2.7.3})$. $(\ref{2.7.2})$ can be proved in a similar fashion. $\Box$\vspace{5mm}

\begin{proposition}\label{p2.3}
Suppose $\hat{u}_{0}$ is supported on $|\xi| \sim N_{1}$ and $\hat{v}_{0}$ is supported on $|\xi| \sim N_{2}$, $N_{1} << N_{2}$. Then

\begin{equation}\label{2.8}
\| (e^{\pm it \Delta} u_{0})(e^{\pm it \Delta} v_{0}) \|_{L_{t,x}^{3}(\mathbf{R} \times \mathbf{R})}
\lesssim (\frac{N_{1}}{N_{2}})^{1/4} \| u_{0} \|_{L^{2}(\mathbf{R})} \| v_{0} \|_{L^{2}(\mathbf{R})}.
\end{equation}
\end{proposition}

\noindent \emph{Proof:} By proposition $\ref{p2.2}$,

$$\| (e^{\pm it \Delta} u_{0})(e^{\pm it \Delta} v_{0}) \|_{L_{t,x}^{2}(\mathbf{R} \times \mathbf{R})}
\lesssim (\frac{1}{N_{2}})^{1/2} \| u_{0} \|_{L^{2}(\mathbf{R})} \| v_{0} \|_{L^{2}(\mathbf{R})}.$$

\noindent Also, combining Strichartz estimates and the Sobolev embedding theorem,

$$\| (e^{\pm it \Delta} u_{0})(e^{\pm it \Delta} v_{0}) \|_{L_{t,x}^{6}(\mathbf{R} \times \mathbf{R})}
\lesssim \| e^{\pm it \Delta} u_{0} \|_{L_{t,x}^{\infty}(\mathbf{R} \times \mathbf{R})} \| e^{\pm it \Delta} u_{0} \|_{L_{t,x}^{6}(\mathbf{R} \times \mathbf{R})} \lesssim N_{1}^{1/2} \| u_{0} \|_{L^{2}(\mathbf{R})} \| v_{0} \|_{L^{2}(\mathbf{R})}. $$

\noindent The proposition follows by interpolation. $\Box$\vspace{5mm}

\noindent Right now, we know that our minimal mass blowup solution is concentrated in both space and frequency, that is,

\begin{equation}\label{9.2}
\int_{|x - x(t)| \geq \frac{C(\eta)}{N(t)}} |u(t,x)|^{2} dx < \eta,
\end{equation}

\begin{equation}\label{9.3}
\int_{|\xi - \xi(t)| \geq C(\eta) N(t)} |\hat{u}(t,\xi)|^{2} d\xi < \eta.
\end{equation}

\noindent Since we will be using the interaction Morawetz estimate, we will not need to track the movement of $x(t)$, however, it will be very important to track the movement of $\xi(t)$. One weapon to partially counter the movement of $\xi(t)$ is the Galilean transformation.

\begin{theorem}\label{t9.2}
Suppose $u(t,x)$ solves

\begin{equation}\label{9.4}
\aligned
i u_{t} + \Delta u &= F(u), \\
u(0,x) &= u_{0}.
\endaligned
\end{equation}

\noindent Then $v(t,x) = e^{-it |\xi_{0}|^{2}} e^{ix \cdot \xi_{0}} u(t, x - 2 \xi_{0} t)$ solves the initial value problem

\begin{equation}\label{9.5}
\aligned
i v_{t} + \Delta v &= F(v), \\
v(0,x) &= e^{ix \cdot \xi_{0}} u(0,x).
\endaligned
\end{equation}

\end{theorem}

\noindent \emph{Proof:} This follows by direct calculation. $\Box$\vspace{5mm}

\noindent If $u(t,x)$ obeys $(\ref{9.2})$ and $(\ref{9.3})$ and $v(t,x) = e^{-it |\xi_{0}|^{2}} e^{ix \cdot \xi_{0}} u(t, x - 2 \xi_{0} t)$, then

\begin{equation}\label{9.7}
\int_{|\xi - \xi_{0} - \xi(t)| \geq C(\eta) N(t)} |\hat{v}(t,\xi)|^{2} d\xi < \eta,
\end{equation}

\begin{equation}\label{9.8}
\int_{|x - 2 \xi_{0} t - x(t)| \geq \frac{C(\eta)}{N(t)}} |v(t,x)|^{2} dx < \eta.
\end{equation}

\noindent \textbf{Remark:} This will be useful to us later because it shifts $\xi(t)$ by a fixed amount $\xi_{0} \in \mathbf{R}^{d}$. For example, this allows us to set $\xi(0) = 0$.

\begin{lemma}\label{l9.1.1}
If $J$ is an interval with

\begin{equation}\label{9.1.1}
\| u \|_{L_{t,x}^{6}(J \times \mathbf{R})} \leq C,
\end{equation}

\noindent then for $t_{1}, t_{2} \in J$,

\begin{equation}\label{9.1.2}
N(t_{1}) \sim_{C, m_{0}} N(t_{2}).
\end{equation}
\end{lemma}

\noindent \emph{Proof:} See \cite{KVZ}, \cite{KTV}, or \cite{TVZ}. $\Box$\vspace{5mm}

\noindent Now if $\| u \|_{L_{t,x}^{6}([0, T] \times \mathbf{R}^{d})} \leq C$, partition $[0, T]$ into $\sim \frac{C^{6}}{\epsilon_{0}^{6}}$ subintervals and iterate. $\Box$ \vspace{5mm}

\noindent We can control the movement of $\xi(t)$ with a similar argument.

\begin{lemma}\label{l9.1.2}
Partition $J = [0, T_{0}]$ into subintervals $J = \cup J_{k}$ such that

\begin{equation}\label{9.1.7}
\| u \|_{L_{t,x}^{6}(J_{k} \times \mathbf{R}^{d})} \leq \epsilon_{0}.
\end{equation}

\noindent Let $N(J_{k}) = \sup_{t \in J_{k}} N(t)$. Then

\begin{equation}\label{9.1.8}
|\xi(0) - \xi(T_{0})| \lesssim \sum_{k} N(J_{k}),
\end{equation}

\noindent which is the sum over the intervals $J_{k}$.
\end{lemma}

\noindent \emph{Proof:} Again take $\eta = \frac{m_{0}^{2}}{1000}$. Let $t_{1}, t_{2} \in J_{k}$. By Strichartz estimates,

\begin{equation}\label{9.1.9}
\| \int_{t_{1}}^{t} e^{i(t - \tau) \Delta} |u(\tau)|^{4} u(\tau) d\tau \|_{L_{x}^{2}(\mathbf{R})} \leq \frac{m_{0}}{1000}.
\end{equation}

\noindent By $(\ref{9.2})$ and $(\ref{9.3})$

\begin{equation}\label{9.1.10}
\int_{|\xi - \xi(t_{1})| \geq C(\frac{m_{0}^{2}}{1000}) N(t_{1})} |\hat{u}(t_{1}, \xi)|^{2} d\xi \leq \frac{m_{0}^{2}}{1000},
\end{equation}

\noindent and

\begin{equation}\label{9.1.11}
\int_{|\xi - \xi(t_{2})| \geq C(\frac{m_{0}^{2}}{1000}) N(t_{2})} |\hat{u}(t_{2}, \xi)|^{2} d\xi \leq \frac{m_{0}^{2}}{1000}.
\end{equation}

\noindent By Duhamel's formula, conservation of mass, $(\ref{9.1.9})$, $(\ref{9.1.10})$, and $(\ref{9.1.11})$, the balls $|\xi - \xi(t)| \leq C(\frac{m_{0}^{2}}{1000}) N(t_{1})$, $|\xi - \xi(t)| \leq C(\frac{m_{0}^{2}}{1000}) N(t_{2})$ must intersect, $|\xi(t_{1}) - \xi(t_{2})| \leq 3 C(\frac{m_{0}^{2}}{1000})(N(t_{1}) + N(t_{2}))$. By the triangle inequality and lemma $\ref{l9.1.1}$,

\begin{equation}\label{9.1.12}
|\xi(T_{0}) - \xi(0)| \leq \sum_{k} |\xi(t_{k}) - \xi(t_{k + 1})| \lesssim \sum_{k} N(t_{k}).
\end{equation}

\noindent $\Box$\vspace{5mm}

\noindent Next, we quote a result,

\begin{lemma}\label{l9.1.2.1}
If $u(t,x)$ is a minimal mass blowup solution on an interval J,

\begin{equation}\label{9.1.12.1}
\int_{J} N(t)^{2} dt \lesssim \| u \|_{L_{t,x}^{6}(J \times \mathbf{R})}^{6} \lesssim 1 + \int_{J} N(t)^{2} dt.
\end{equation}
\end{lemma}

\noindent \emph{Proof:} See \cite{KVZ}.\vspace{5mm}

\noindent Finally we will prove a lemma that will be useful to us when analyzing the blowup scenarios with $N(t) \leq 1$.

\begin{lemma}\label{l9.1.4}
Suppose $u$ is a minimal mass blowup solution with $N(t) \leq 1$. Suppose also that $J$ is some interval partitioned into subintervals $J_{k}$ with $\| u \|_{L_{t,x}^{6}(J_{k} \times \mathbf{R})} = \epsilon_{0}$ on each $J_{k}$. Again let

\begin{equation}\label{9.1.17.1}
N(J_{k}) = \sup_{J_{k}} N(t).
\end{equation}

\noindent Then,

\begin{equation}\label{9.1.18}
\sum_{J_{k}} N(J_{k}) \sim \int_{J} N(t)^{3} dt.
\end{equation}
\end{lemma}

\noindent \emph{Proof:} By lemma $\ref{l9.1.2.1}$,

\begin{equation}\label{9.1.19}
\int_{J} N(t)^{2} \lesssim \| u \|_{L_{t,x}^{6}(J \times \mathbf{R})}^{6}.
\end{equation}

\noindent Since $\| u \|_{L_{t,x}^{6}(J_{k} \times \mathbf{R})} = \epsilon_{0}$, by $(\ref{9.1.12.1})$,

$$\int_{J_{k}} N(t)^{3} dt \lesssim N(J_{k}) \int_{J_{k}} N(t)^{2} \lesssim \epsilon_{0}^{6} N(J_{k}),$$

\noindent so $$\int_{J} N(t)^{3} dt \lesssim \sum_{J_{k}} N(J_{k}).$$

\noindent On the other hand, by the Duhamel formula,

\begin{equation}\label{9.1.19.2}
\| u \|_{L_{t}^{4} L_{x}^{\infty}(J_{k} \times \mathbf{R})} \lesssim \| u_{0} \|_{L^{2}(\mathbf{R})} + \| u \|_{L_{t,x}^{6}(J_{k} \times \mathbf{R})}^{5} \lesssim 1.
\end{equation}

\noindent Interpolating this with

\begin{equation}\label{9.1.19.3}
\| u_{|\xi - \xi(t)| \geq C(\eta) N(t)} \|_{L_{t}^{\infty} L_{x}^{2}(J_{k} \times \mathbf{R})} \leq \eta^{1/2},
\end{equation}

\noindent we have

\begin{equation}\label{9.1.19.4}
\| u_{|\xi - \xi(t)| \geq C(\eta(\epsilon)) N(t)} \|_{L_{t,x}^{6}(J_{k} \times \mathbf{R})} \leq \frac{\epsilon_{0}}{1000},
\end{equation}

\noindent for a small, fixed $\eta(\epsilon_{0}) > 0$. By the Sobolev embedding theorem,

\begin{equation}\label{9.1.19.5}
\| u_{|\xi - \xi(t)| \leq C(\eta(\epsilon_{0})) N(t)}(t) \|_{L_{x}^{6}(\mathbf{R}^{2})} \lesssim [C(\eta(\epsilon_{0})) N(t)]^{\frac{1}{3}}.
\end{equation}

\noindent Therefore, $$\epsilon_{0}^{6} \lesssim \int_{J_{k}} C(\eta(\epsilon_{0}))^{2} N(t)^{2} dt.$$ Since $N(t_{1}) \sim N(t_{2})$ for $t_{1}, t_{2} \in J_{k}$, this implies

\begin{equation}\label{9.1.19.6}
N(J_{k}) \lesssim \int_{J_{k}} N(t)^{3} dt.
\end{equation}

\noindent Summing up over subintervals proves the lemma. $\Box$

\section{A norm adapted to $\xi(t)$, $N(t)$ constant}
\noindent As a warm-up, we will treat the minimal mass blowup scenario $N(t) \equiv 1$, $\xi(t) \equiv 0$, $\mu = \mp 1$. Rescaling,

\begin{equation}\label{3.1}
u_{\lambda}(t,x) = \frac{1}{\lambda^{1/2}} u(\frac{t}{\lambda^{2}}, \frac{x}{\lambda}),
\end{equation}

\noindent $N(t) \equiv \frac{1}{\lambda}$. We will choose to treat the case $N(t) = \delta$, $\delta > 0$ sufficiently small so that $\delta < \epsilon^{10}$, and dropping the $\lambda$ from $u_{\lambda}(t)$,

\begin{equation}\label{3.2}
\| P_{> \frac{\delta^{1/2}}{32}} u(t) \|_{L_{t}^{\infty} L_{x}^{2}((-\infty, \infty) \times \mathbf{R}^{2})} < \epsilon,
\end{equation}

\noindent and for any $a \in \mathbf{R}$,

\begin{equation}\label{3.3}
\| u \|_{L_{t}^{4} L_{x}^{\infty}([a, a + 1] \times \mathbf{R}^{2})} + \| u \|_{L_{t,x}^{6}([a, a + 1] \times \mathbf{R}^{2})} \leq \epsilon_{0}.
\end{equation}

\noindent The semi-norm we are about to define is adapted to the case $N(t) \equiv \delta$. This semi-norm will be generalized in the next section to treat the case when $N(t)$ and $\xi(t)$ are free to move around.

\begin{definition}\label{d3.1}
\noindent Let $N_{j}$ be a dyadic integer.

\begin{equation}\label{3.6}
\aligned
\| u \|_{X_{N_{j}}^{k}}^{2} = \sum_{1 \leq N_{i} \leq N_{j}} \frac{N_{i}}{N_{j}} \sum_{l = 0}^{\frac{N_{j}}{N_{i}} - 1} \| P_{N_{i}} u \|_{U_{\Delta}^{2}([k N_{j} + l N_{i}, k N_{j} + (l + 1) N_{i}] \times \mathbf{R})}^{2} \\
+ \sum_{N_{j} < N_{i}} \| P_{N_{i}} u \|_{U_{\Delta}^{2}([k N_{j}, (k + 1) N_{j}] \times \mathbf{R})}^{2}.
\endaligned
\end{equation}

\noindent Now let $M$ be some dyadic integer,

\begin{equation}\label{3.7}
\| u \|_{X_{M}([0, M] \times \mathbf{R})}^{2} \equiv \sup_{1 \leq N_{j} \leq M} \sup_{0 \leq k \leq \frac{M}{N_{j}}} \| u \|_{X_{N_{j}}^{k}}^{2}.
\end{equation}

\noindent Similarly, we can define

\begin{equation}\label{3.8}
\| u \|_{X_{M}([a, a + M] \times \mathbf{R})}^{2}
\end{equation}

\noindent for any $a \in \mathbf{R}$.
\end{definition}

\noindent \textbf{Remark:} $\| u \|_{X_{M}([0, M] \times \mathbf{R})}$ is only a semi-norm since if $f(t)$ is a nonzero function supported on $|\xi| < 1$, $\| f(t) \|_{X_{M}([0, M] \times \mathbf{R})} \equiv 0$. Therefore, we need to say something about a minimal mass blowup solution at low frequencies.

\begin{lemma}\label{l3.1}
Suppose $u(t)$ is a minimal mass blowup solution to $(1.1)$, $\mu = \pm 1$, and $J$ is an interval with

\begin{equation}\label{3.8.1}
\| u \|_{L_{t,x}^{6}(J \times \mathbf{R})} + \| u \|_{L_{t}^{4} L_{x}^{\infty}(J \times \mathbf{R})} \leq \epsilon_{0},
\end{equation}

\noindent and $N(t) = \delta$ on $J$. Then

\begin{equation}\label{3.8.2}
\| P_{> \delta^{1/2}} u(t) \|_{U_{\Delta}^{2}(J \times \mathbf{R})} \lesssim \epsilon.
\end{equation}
\end{lemma}

\noindent \emph{Proof:} Let $J = [a, b]$. By Duhamel's formula, for $t \in J$,

\begin{equation}\label{3.8.3}
u(t) = e^{i(t - a) \Delta} u(a) - i \int_{a}^{t} e^{i(t - \tau) \Delta} F(u(\tau)) d\tau.
\end{equation}

\noindent Since $$\| P_{> \frac{\delta^{1/2}}{32}} u(t) \|_{L_{t}^{\infty} L_{x}^{2}(J \times \mathbf{R})} \leq \epsilon,$$

$$\| P_{> \frac{\delta^{1/2}}{32}} e^{i(t - a) \Delta} u(a) \|_{U_{\Delta}^{2}(J \times \mathbf{R})} \lesssim \epsilon.$$

\noindent Also,

$$\| \int_{a}^{t} e^{i(t - \tau) \Delta} P_{> \delta^{1/2}}(|u(\tau)|^{4} u(\tau)) d\tau \|_{U_{\Delta}^{2}(J \times \mathbf{R})} \lesssim \| P_{> \delta^{1/2}} (|u(\tau)|^{4} u(\tau)) \|_{L_{t}^{4/3} L_{x}^{1}(J \times \mathbf{R})}$$

$$\lesssim \| P_{> \frac{\delta^{1/2}}{32}} u \|_{L_{t}^{\infty} L_{x}^{2}(J \times \mathbf{R})} \| u \|_{L_{t}^{16/3} L_{x}^{4}(J \times \mathbf{R})}^{4} \lesssim \epsilon_{0}^{4} \epsilon. \Box$$

\noindent \textbf{Remark:} Using the exact same arguments, if $J$ is an interval with

\begin{equation}\label{3.8.4}
\| u \|_{L_{t,x}^{6}(J \times \mathbf{R})} + \| u \|_{L_{t}^{4} L_{x}^{\infty}(J \times \mathbf{R})} \leq \epsilon_{0},
\end{equation}

\begin{equation}\label{3.8.5}
\| P_{> N(J) \delta^{1/2}} u \|_{U_{\Delta}^{2}(J \times \mathbf{R})} \lesssim \epsilon.
\end{equation}

\begin{theorem}\label{t3.1}
Suppose $u(t)$ is a minimal mass blowup solution to

\begin{equation}\label{3.9}
i u_{t} + \Delta u = F(u).
\end{equation}

\noindent There exists a fixed constant $C$ such that for $\epsilon, \delta(\epsilon) > 0$ sufficiently small,

\begin{equation}\label{3.10}
\| u \|_{X_{M}([0, M] \times \mathbf{R})} \leq C \epsilon
\end{equation}

\noindent for all dyadic $M$, $1 \leq M < \infty$.
\end{theorem}

\noindent \emph{Sketch of Proof:} Theorem $\ref{t3.1}$ is proved by induction. By lemma $\ref{l3.1}$,

\begin{equation}\label{3.10.1}
\| P_{1} u(t) \|_{U_{\Delta}^{2}([a, a + 1] \times \mathbf{R})} \leq C \epsilon.
\end{equation}

\noindent Suppose that for any dyadic integer $M$, $1 \leq M < \infty$,

\begin{equation}\label{3.11}
\| u \|_{X_{M}([a, a + M] \times \mathbf{R})} \leq \frac{C}{2} \epsilon + \frac{C}{2}(\epsilon^{2} + \| u \|_{X_{M}([a, a + M] \times \mathbf{R}^{2})}^{2}),
\end{equation}

\noindent $C$ is independent of $\epsilon > 0$. Then we are able to prove theorem $\ref{t3.1}$ by induction. Suppose that for $M \leq N$,

\begin{equation}\label{3.12}
\| u(t) \|_{X_{M}([a, a + M] \times \mathbf{R})} \leq C \epsilon,
\end{equation}

\noindent for a fixed constant $C$, $\epsilon > 0$, and for any $a \in \mathbf{R}$. Then making a crude estimate,

\begin{equation}\label{3.13}
\| u(t) \|_{X_{2N}([a, a + 2N] \times \mathbf{R})} \leq 2 C \epsilon.
\end{equation}

\noindent By $(\ref{3.11})$, $(\ref{3.13})$,

\begin{equation}\label{3.14}
\| u(t) \|_{X_{2N}([a, a + 2N] \times \mathbf{R})} \leq \frac{C}{2} \epsilon + \frac{C}{2}(\epsilon^{2} + (2C \epsilon)^{2}).
\end{equation}

\noindent For $\epsilon > 0$ sufficiently small, this implies that for $a \in \mathbf{R}$,

\begin{equation}\label{3.15}
\| u(t) \|_{X_{2N}([a, a + 2N] \times \mathbf{R})} \leq C \epsilon,
\end{equation}

\noindent closing the induction. $\Box$\vspace{5mm}

\noindent The proof of an estimate of the form $(\ref{3.11})$ will occupy the bulk of the paper. In fact, we will prove an estimate of the form $(\ref{3.11})$ for a generalization of $\| u \|_{X_{M}([a, a + M] \times \mathbf{R})}$ used to treat the case when $N(t)$ need not be constant. $(\ref{3.11})$ will be a special case of the more general result.\vspace{5mm}

\noindent The purpose of this section is to discuss the simpler case in the hopes that the main idea is more evident, since it is not obscured by the technical details that arise when $\xi(t)$ and $N(t)$ are free to move around.

\section{Estimates when $N(t)$, $\xi(t)$ are free to vary}
In this section we will generalize the seminorm in the previous section to adapt it to the case when $N(t)$ and $\xi(t)$ are free to vary. We will define the seminorm $\tilde{X}_{M}([0, T])$ on the time interval $[0, T]$ to be an analogue of the $X_{M}([0, M])$ norm defined in the previous section.\vspace{5mm}

\noindent Suppose $[0, T] = \cup_{l = 1}^{M} J_{l}$, with $\| u \|_{L_{t,x}^{6}(J_{l} \times \mathbf{R})} = \epsilon_{0}$, $\sum_{J_{l}} N(J_{l}) = \delta M$. We will call the individual $J_{l}$ subintervals the small intervals. We want to partition $[0, T]$ at level $N_{i}$ for $1 \leq N_{i} \leq M$. If $N(J_{l}) > \frac{\delta N_{i}}{2}$ then we will call $J_{l}$ a red interval at level $N_{i}$.\vspace{5mm}

\noindent A union $G = \cup J_{l}$ of $N_{i}$ consecutive small intervals with $$\sum_{J_{l} \subset G} N(J_{l}) \leq \delta N_{i}$$ and $N(J_{l}) \leq \frac{\delta N_{i}}{2}$ for each $J_{l} \subset G$ will be called a length green interval at level $N_{i}$. A union $G$ of $\leq N_{i}$ consecutive small intervals $J_{l}$ with $$\frac{\delta N_{i}}{2} < \sum_{J_{l} \subset G} N(J_{l}) \leq \delta N_{i}$$ will be called a weight green interval at level $N_{i}$.\vspace{5mm}

\noindent A union $Y$ of $< N_{i}$ consecutive small intervals with $$\sum_{J_{l} \subset Y} N(J_{l}) \leq \frac{\delta N_{i}}{2}$$ will be called a yellow interval at level $N_{i}$.\vspace{5mm}

\noindent $[0, T]$ will be partitioned so that every yellow interval $Y$ lies immediately to the left of a red interval, or that $T \in Y$. If there is a yellow interval $Y$, and the small interval $J_{l}$ to the right of $Y$ satisfies $N(J_{l}) \leq \frac{\delta N_{i}}{2}$ then we take $Y \cup J_{l} = Y^{\ast}$. $Y^{\ast}$ is the union of $\leq N_{i}$ small intervals with $$\sum_{J_{l} \subset Y^{\ast}} N(J_{l}) \leq \delta N_{i}.$$ If $Y^{\ast}$ is the union of $N_{i}$ small intervals then $Y^{\ast}$ is a length green interval. If $$\frac{\delta N_{i}}{2} < \sum_{J_{l} \subset Y^{\ast}} N(J_{l}) \leq \delta N_{i},$$ then $Y^{\ast}$ is a weight green interval. If $$\sum_{J_{l} \subset Y^{\ast}} N(J_{l}) \leq \frac{\delta N_{i}}{2}$$ and $Y^{\ast}$ is the union of $< N_{i}$ small intervals, then $Y^{\ast}$ remains a yellow interval. If $T \notin Y^{\ast}$ and the small interval to the right of $Y^{\ast}$ is not red, then repeat the above procedure.\vspace{5mm}

\noindent \textbf{Remark:} We will always say that $[0, T]$ is a green interval at level $M$.\vspace{5mm}

\noindent \textbf{Remark:} The reader should think of the yellow intervals at level $N_{i}$ as the scraps left over after carving out the red and green intervals at level $N_{i}$.\vspace{5mm}

\noindent We also want to apply the seminorms in the previous section to the case when $\xi(t)$ is free to travel around in $\mathbf{R}$. This seminorm was defined in \cite{D3} for any dimension $d$, $d \geq 1$.

\begin{definition}\label{d5.1}
\noindent For a green interval $G_{\alpha}^{N_{i}} = [a, b]$, let $\xi(G_{\alpha}^{N_{i}}) = \xi(a)$. $\xi(Y_{\alpha}^{N_{i}})$ and $\xi(R_{\alpha}^{N_{i}})$ can be defined in a similar manner.\vspace{5mm}

\noindent If $G_{k}^{N_{j}}$ is a green interval at level $N_{j}$, then

\begin{equation}\label{5.1}
\aligned
\| u \|_{X(G_{k}^{N_{j}})}^{2} \equiv \sum_{1 \leq N_{i} \leq N_{j}} (\frac{N_{i}}{N_{j}}) \sum_{G_{\alpha}^{N_{i}} \cap G_{k}^{N_{j}} \neq \emptyset} \| \tilde{P}_{\xi(G_{\alpha}^{N_{i}}), \frac{N_{i}}{4} \leq \cdot \leq 4 N_{i}} u \|_{U_{\Delta}^{2}(G_{\alpha}^{N_{i}} \times \mathbf{R}^{d})}^{2} \\
+ \sum_{N_{j} < N_{i}} \| \tilde{P}_{\xi(G_{k}^{N_{j}}), \frac{N_{i}}{4} \leq \cdot \leq 4 N_{i}} u \|_{U_{\Delta}^{2}(G_{k}^{N_{j}} \times \mathbf{R}^{d})}^{2}
+ \sup_{1 \leq N_{i} < N_{j}} \sup_{Y_{\alpha}^{N_{i}} \cap G_{k}^{N_{j}} \neq \emptyset} \| \tilde{P}_{\xi(Y_{\alpha}^{N_{i}}), \frac{N_{i}}{4} \leq \cdot \leq 4 N_{i}} u \|_{U_{\Delta}^{2}(Y_{\alpha}^{N_{i}} \times \mathbf{R}^{d})}^{2}.
\endaligned
\end{equation}

\noindent For a yellow interval at level $N_{j}$,

\begin{equation}\label{5.2}
\aligned
\| u \|_{X(Y_{k}^{N_{j}})}^{2} \equiv \sum_{1 \leq N_{i} \leq N_{j}} (\frac{N_{i}}{N_{j}}) \sum_{G_{\alpha}^{N_{i}} \cap Y_{k}^{N_{j}} \neq \emptyset} \| P_{\xi(G_{\alpha}^{N_{i}}), \frac{N_{i}}{4} \leq \cdot \leq 4 N_{i}} u \|_{U_{\Delta}^{2}(G_{\alpha}^{N_{i}} \times \mathbf{R}^{d})}^{2} \\
+ \sum_{N_{j} \leq N_{i}} \| P_{\xi(Y_{k}^{N_{j}}), \frac{N_{i}}{4} \leq \cdot \leq 4 N_{i}} u \|_{U_{\Delta}^{2}(Y_{k}^{N_{j}} \times \mathbf{R}^{d})}^{2}
+ \sup_{1 \leq N_{i} < N_{j}} \sup_{Y_{\alpha}^{N_{i}} \cap Y_{k}^{N_{j}} \neq \emptyset} \| P_{\xi(Y_{\alpha}^{N_{i}}), \frac{N_{i}}{4} \leq \cdot \leq 4 N_{i}} u \|_{U_{\Delta}^{2}(Y_{\alpha}^{N_{i}} \times \mathbf{R}^{d})}^{2}.
\endaligned
\end{equation}

\noindent Then,

\begin{equation}\label{5.3}
\| u \|_{\tilde{X}_{M}([0, T] \times \mathbf{R}^{d})}^{2} \equiv \sup_{1 \leq N_{j} \leq M} \sup_{G_{k}^{N_{j}} \subset [0, T]} \| u \|_{X(G_{k}^{N_{j}})}^{2} + \sup_{1 \leq N_{j} \leq M} \sup_{Y_{k}^{N_{j}} \subset [0, T]} \| u \|_{X(Y_{k}^{N_{j}})}^{2}.
\end{equation}

\noindent Also for a dyadic integer $N$, $1 \leq N \leq M$ define the norm

\begin{equation}\label{5.4}
\| u \|_{\tilde{X}_{N}([0, T] \times \mathbf{R}^{d})}^{2} \equiv \sup_{1 \leq N_{j} \leq N} \sup_{G_{k}^{N_{j}} \subset [0, T]} \| u \|_{X(G_{k}^{N_{j}})}^{2} + \sup_{1 \leq N_{j} \leq N} \sup_{Y_{k}^{N_{j}} \subset [0, T]} \| u \|_{X(Y_{k}^{N_{j}})}^{2}.
\end{equation}
\end{definition}

\noindent We first prove than an estimate on $\| u \|_{X(G_{k}^{N_{j}})}$ gives control over $\| P_{\xi(G_{k}^{N_{j}}, N_{i}} u \|_{U_{\Delta}^{2}(G_{k}^{N_{j}} \times \mathbf{R})}$ for a dyadic frequency $N_{i}$, along with Strichartz estimates of $P_{\xi(G_{k}^{N_{j}}), N_{i}} u$.

\begin{lemma}\label{l5.0}
For a dyadic frequency $1 \leq N_{j}$,

\begin{equation}\label{5.5}
\aligned
\sum_{1 \leq N_{i} \leq N_{j}} (\frac{N_{i}}{N_{j}}) [\sum_{G_{\alpha}^{N_{i}} \cap G_{k}^{N_{j}} \neq \emptyset} \| P_{\xi(G_{\alpha}^{N_{i}}), N_{i}} u \|_{U_{\Delta}^{2}(G_{\alpha}^{N_{i}} \times \mathbf{R}^{d})}^{2} + \sum_{Y_{\alpha'}^{N_{i}} \cap G_{k}^{N_{j}} \neq \emptyset} \| P_{\xi(Y_{\alpha'}^{N_{i}}), N_{i}} u \|_{U_{\Delta}^{2}(Y_{\alpha'}^{N_{i}} \times \mathbf{R}^{d})}^{2}\\
+ \sum_{R_{\alpha''}^{N_{i}} \subset G_{k}^{N_{j}}} \| P_{\xi(R_{\alpha''}^{N_{i}}), N_{i}} u \|_{U_{\Delta}^{2}(R_{\alpha''}^{N_{i}} \times \mathbf{R}^{d})}^{2}] + \sum_{N_{j} < N_{i}} \| P_{\xi(G_{k}^{N_{j}}), N_{i}} u \|_{U_{\Delta}^{2}(G_{k}^{N_{j}} \times \mathbf{R}^{d})}^{2} \lesssim \| u \|_{X(G_{k}^{N_{j}})}^{2} + \epsilon^{2}.
\endaligned
\end{equation}

\noindent Similarly,

\begin{equation}\label{5.6}
\aligned
\sum_{1 \leq N_{i} \leq N_{j}} (\frac{N_{i}}{N_{j}}) [\sum_{G_{\alpha}^{N_{i}} \cap Y_{k}^{N_{j}} \neq \emptyset} \| P_{\xi(G_{\alpha}^{N_{i}}), N_{i}} u \|_{U_{\Delta}^{2}(G_{\alpha}^{N_{i}} \times \mathbf{R}^{d})}^{2} + \sum_{Y_{\alpha'}^{N_{i}} \cap Y_{k}^{N_{j}} \neq \emptyset} \| P_{\xi(Y_{\alpha'}^{N_{i}}), N_{i}} u \|_{U_{\Delta}^{2}(Y_{\alpha'}^{N_{i}} \times \mathbf{R}^{d})}^{2}\\
+ \sum_{R_{\alpha''}^{N_{i}} \subset Y_{k}^{N_{j}}} \| P_{\xi(R_{\alpha''}^{N_{i}}), N_{i}} u \|_{U_{\Delta}^{2}(R_{\alpha''}^{N_{i}} \times \mathbf{R}^{d})}^{2}] + \sum_{N_{j} < N_{i}} \| P_{\xi(Y_{k}^{N_{j}}), N_{i}} u \|_{U_{\Delta}^{2}(Y_{k}^{N_{j}} \times \mathbf{R}^{d})}^{2} \lesssim \| u \|_{X(Y_{k}^{N_{j}})}^{2} + \epsilon^{2}.
\endaligned
\end{equation}

\noindent Finally, for $1 \leq N_{i} < N_{j}$, suppose $(p, q)$ is a $d$-admissible pair.

\begin{equation}\label{5.6.1}
\| P_{\xi(t), N_{i}} u \|_{L_{t}^{p} L_{x}^{q}(G_{k}^{N_{j}} \times \mathbf{R}^{d})} \lesssim (\frac{N_{j}}{N_{i}})^{1/p} (\delta^{1/2p} + \epsilon + \| u \|_{\tilde{X}_{N_{j}}}),
\end{equation}

\noindent and

\begin{equation}\label{5.6.1}
\| P_{\xi(t), N_{i}} u \|_{L_{t}^{p} L_{x}^{q}(Y_{k}^{N_{j}} \times \mathbf{R}^{d})} \lesssim (\frac{N_{j}}{N_{i}})^{1/p} (\delta^{1/2p} + \epsilon + \| u \|_{\tilde{X}_{N_{j}}}).
\end{equation}
\end{lemma}

\noindent \emph{Proof:} See \cite{D3}.\vspace{5mm}

\noindent Also, recall that from \cite{D3}

\begin{theorem}\label{t5.1}
If $u(t)$ is a minimal mass blowup solution to

\begin{equation}\label{5.11}
\aligned
i u_{t} + \Delta u = F(u), \\
u(0,x) = u_{0} \in L^{2}(\mathbf{R}),
\endaligned
\end{equation}

\noindent $\mu = \pm 1$. There exists a constant $C$ such that for $\epsilon > 0$, $\delta(\epsilon) > 0$ sufficiently small, for any dyadic integer $M$, if there exist small intervals $J_{l}$ with $$[0, T] = \cup_{l = 1}^{M} J_{l},$$ $$\sum_{J_{l} \subset [0, T]} N(J_{l}) = \frac{\delta M}{2},$$ $$\| u \|_{L_{t,x}^{6}(J_{l} \times \mathbf{R})} = \epsilon_{0},$$ then

\begin{equation}\label{5.12}
\| u \|_{\tilde{X}_{M}([0, T] \times \mathbf{R})} \leq C \epsilon.
\end{equation}
\end{theorem}

\noindent It was showed in \cite{D3} that to prove theorem $\ref{t5.1}$ it suffices to prove two intermediate lemmas. Lemmas $\ref{l5.2}$ and $\ref{l5.3}$ will be proved in this section and the next.

\begin{lemma}\label{l5.2}
If $u(t)$ satisfies the conditions in theorem $\ref{t5.1}$, then

\begin{equation}\label{5.13}
\| u \|_{X(G_{k}^{N_{j}})}^{2} \lesssim \epsilon^{2}
\end{equation}

\begin{equation}\label{5.14}
+ (\epsilon + \| u \|_{\tilde{X}_{N_{j}}([0, T] \times \mathbf{R})})^{2} \sum_{1 \leq N_{i} \leq N_{j}} (\frac{N_{i}}{N_{j}}) \sum_{G_{\alpha}^{N_{i}} \cap G_{k}^{N_{j}} \neq \emptyset} \sum_{\frac{N_{i}}{32} \leq N_{1} \leq 32 N_{i}} \| P_{\xi(G_{\alpha}^{N_{i}}), N_{1}} u \|_{U_{\Delta}^{2}(G_{\alpha}^{N_{i}} \times \mathbf{R})}^{2}
\end{equation}

\begin{equation}\label{5.15}
+ (\epsilon + \| u \|_{\tilde{X}_{N_{j}}([0, T] \times \mathbf{R})})^{2} \sum_{N_{j} < N_{i}} \sum_{\frac{N_{i}}{32} \leq N_{1} \leq 32 N_{i}} \| P_{\xi(G_{k}^{N_{j}}), N_{1}} u \|_{U_{\Delta}^{2}(G_{k}^{N_{j}} \times \mathbf{R})}^{2}
\end{equation}

\begin{equation}\label{5.16}
+ (\epsilon + \| u \|_{\tilde{X}_{N_{j}}([0, T] \times \mathbf{R})})^{2} \sup_{1 \leq N_{i} < N_{j}} \sup_{Y_{\alpha}^{N_{i}} \cap G_{k}^{N_{j}} \neq \emptyset} \sum_{\frac{N_{i}}{32} \leq N_{1} \leq 32 N_{i}} \| P_{\xi(Y_{\alpha}^{N_{i}}), N_{1}} u \|_{U_{\Delta}^{2}(Y_{\alpha}^{N_{i}} \times \mathbf{R})}^{2}
\end{equation}

\begin{equation}\label{5.17}
+ (\epsilon + \| u \|_{\tilde{X}_{N_{j}}([0, T] \times \mathbf{R})})^{2} \sum_{1 \leq N_{i} \leq N_{j}} (\frac{N_{i}}{N_{j}}) \sum_{G_{\alpha}^{N_{i}} \cap G_{k}^{N_{j}} \neq \emptyset} (\sum_{32 N_{i} < N_{1}} (\frac{N_{i}}{N_{1}})^{1/4} \| P_{\xi(G_{\alpha}^{N_{i}}), N_{1}} u \|_{U_{\Delta}^{2}(G_{\alpha}^{N_{i}} \times \mathbf{R})})^{2}
\end{equation}

\begin{equation}\label{5.18}
+ (\epsilon + \| u \|_{\tilde{X}_{N_{j}}([0, T] \times \mathbf{R})})^{2} \sum_{N_{j} < N_{i}} (\sum_{32 N_{i} < N_{1}} (\frac{N_{i}}{N_{1}})^{1/4} \| P_{\xi(G_{k}^{N_{j}}), N_{1}} u \|_{U_{\Delta}^{2}(G_{k}^{N_{j}} \times \mathbf{R})})^{2}
\end{equation}

\begin{equation}\label{5.19}
+ (\epsilon + \| u \|_{\tilde{X}_{N_{j}}([0, T] \times \mathbf{R})})^{2} \sup_{1 \leq N_{i} < N_{j}} \sup_{Y_{\alpha}^{N_{i}} \cap G_{k}^{N_{j}} \neq \emptyset} (\sum_{32 N_{i} < N_{1}} (\frac{N_{i}}{N_{1}})^{1/4} \| P_{\xi(Y_{\alpha}^{N_{i}}), N_{1}} u \|_{U_{\Delta}^{2}(Y_{\alpha}^{N_{i}} \times \mathbf{R})})^{2}.
\end{equation}

\end{lemma}

\noindent Similarly, we prove

\begin{lemma}\label{l5.3}
If $u(t)$ satisfies the conditions in theorem $\ref{t5.1}$, then

\begin{equation}\label{5.20}
\| u \|_{X(Y_{k}^{N_{j}})}^{2} \lesssim \epsilon^{2}
\end{equation}

\begin{equation}\label{5.21}
+ (\epsilon + \| u \|_{\tilde{X}_{N_{j}}([0, T] \times \mathbf{R})})^{2} \sum_{1 \leq N_{i} \leq N_{j}} (\frac{N_{i}}{N_{j}}) \sum_{G_{\alpha}^{N_{i}} \cap Y_{k}^{N_{j}} \neq \emptyset} \sum_{\frac{N_{i}}{32} \leq N_{1} \leq 32 N_{i}} \| P_{\xi(G_{\alpha}^{N_{i}}), N_{1}} u \|_{U_{\Delta}^{2}(G_{\alpha}^{N_{i}} \times \mathbf{R})}^{2}
\end{equation}

\begin{equation}\label{5.22}
+ (\epsilon + \| u \|_{\tilde{X}_{N_{j}}([0, T] \times \mathbf{R})})^{2} \sum_{N_{j} < N_{i}} \sum_{\frac{N_{i}}{32} \leq N_{1} \leq 32 N_{i}} \| P_{\xi(Y_{k}^{N_{j}}), N_{1}} u \|_{U_{\Delta}^{2}(Y_{k}^{N_{j}} \times \mathbf{R})}^{2}
\end{equation}

\begin{equation}\label{5.23}
+ (\epsilon + \| u \|_{\tilde{X}_{N_{j}}([0, T] \times \mathbf{R})})^{2} \sup_{1 \leq N_{i} < N_{j}} \sup_{Y_{\alpha}^{N_{i}} \cap Y_{k}^{N_{j}} \neq \emptyset} \sum_{\frac{N_{i}}{32} \leq N_{1} \leq 32 N_{i}} \| P_{\xi(Y_{\alpha}^{N_{i}}), N_{1}} u \|_{U_{\Delta}^{2}(Y_{\alpha}^{N_{i}} \times \mathbf{R})}^{2}
\end{equation}

\begin{equation}\label{5.24}
+ (\epsilon + \| u \|_{\tilde{X}_{N_{j}}([0, T] \times \mathbf{R})})^{2} \sum_{1 \leq N_{i} \leq N_{j}} (\frac{N_{i}}{N_{j}}) \sum_{G_{\alpha}^{N_{i}} \cap Y_{k}^{N_{j}} \neq \emptyset} (\sum_{32 N_{i} < N_{1}} (\frac{N_{i}}{N_{1}})^{1/4} \| P_{\xi(G_{\alpha}^{N_{i}}), N_{1}} u \|_{U_{\Delta}^{2}(G_{\alpha}^{N_{i}} \times \mathbf{R})})^{2}
\end{equation}

\begin{equation}\label{5.25}
+ (\epsilon + \| u \|_{\tilde{X}_{N_{j}}([0, T] \times \mathbf{R})})^{2} \sum_{N_{j} < N_{i}} (\sum_{32 N_{i} < N_{1}} (\frac{N_{i}}{N_{1}})^{1/4} \| P_{\xi(Y_{k}^{N_{j}}), N_{1}} u \|_{U_{\Delta}^{2}(Y_{k}^{N_{j}} \times \mathbf{R})})^{2}
\end{equation}

\begin{equation}\label{5.26}
+ (\epsilon + \| u \|_{\tilde{X}_{N_{j}}([0, T] \times \mathbf{R})})^{2} \sup_{1 \leq N_{i} \leq N_{j}} \sup_{Y_{\alpha}^{N_{i}} \cap Y_{k}^{N_{j}} \neq \emptyset} (\sum_{32 N_{i} < N_{1}} (\frac{N_{i}}{N_{1}})^{1/4} \| P_{\xi(Y_{\alpha}^{N_{i}}), N_{1}} u \|_{U_{\Delta}^{2}(Y_{\alpha}^{N_{i}} \times \mathbf{R})})^{2}.
\end{equation}

\end{lemma}

\noindent \emph{Start of the proof of lemma $\ref{l5.2}$ and lemma $\ref{l5.3}$:} Take a yellow interval $Y_{\alpha'}^{N_{i}}$. For any $a_{\alpha', N_{i}} \in Y_{\alpha'}^{N_{i}}$, the solution to $(\ref{0.1})$ on $Y_{\alpha'}^{N_{i}}$ is equal to

\begin{equation}\label{5.101}
e^{i(t - a_{\alpha', N_{i}}) \Delta} u(a_{\alpha', N_{i}}) - i \int_{a_{\alpha'}^{N_{i}}}^{t} e^{i(t - \tau) \Delta} F(u(\tau)) d\tau.
\end{equation}

\noindent We will postpone the treatment of the Duhamel term, $$\int_{a_{\alpha'}^{N_{i}}}^{t} e^{i(t - \tau) \Delta} |u(\tau)|^{4} u(\tau) d\tau,$$ until the next section. Since $N(t) \leq \frac{\delta N_{i}}{2}$ on $Y_{\alpha}^{N_{i}}$, $$\| P_{\xi(Y_{\alpha}^{N_{i}}), \frac{N_{i}}{4} \leq \cdot \leq 4 N_{i}} u(a_{\alpha, N_{i}}) \|_{L_{x}^{2}(\mathbf{R})} \leq \epsilon.$$

\noindent Next take $G_{L}^{N_{i}}$ and $G_{R}^{N_{i}}$. For $a_{L, N_{i}} \in G_{L}^{N_{i}}$ and $a_{R, N_{i}} \in G_{R}^{N_{i}}$,

$$\sum_{1 \leq N_{i} < N_{j}} (\frac{N_{i}}{N_{j}}) (\| P_{\xi(G_{L}^{N_{i}}, \frac{N_{i}}{4} \leq \cdot \leq 4 N_{i}} u(a_{L, N_{i}}) \|_{L_{x}^{2}(\mathbf{R})}^{2} + \| P_{G_{R}^{N_{i}}, \frac{N_{i}}{4} \leq \cdot \leq 4 N_{i}} u(a_{R, N_{i}}) \|_{L_{x}^{2}(\mathbf{R})}^{2}) \lesssim \epsilon^{2}.$$

\noindent Finally, if $G_{\alpha}^{N_{i}} \subset G_{k}^{N_{j}}$ is a length green interval, we can choose $a_{\alpha, N_{i}}$ such that

\begin{equation}\label{5.102}
\| P_{\xi(G_{\alpha}^{N_{i}}), \frac{N_{i}}{4} \leq \cdot \leq 4 N_{i}} u(a_{\alpha, N_{i}}) \|_{L_{x}^{2}(\mathbf{R})}^{2} \leq \frac{1}{N_{i}} \sum_{J_{l} \subset G_{\alpha}^{N_{i}}} \| P_{\xi(a_{l}), \frac{N_{i}}{8} \leq \cdot \leq 8 N_{i}} u(a_{l}) \|_{L_{x}^{2}(\mathbf{R})}^{2},
\end{equation}

\noindent where $J_{l} = [a_{l}, b_{l}]$. If $G_{\alpha, N_{i}}$ is a weight green interval then we can choose $a_{\alpha, N_{i}} \in G_{\alpha, N_{i}}$ such that

\begin{equation}\label{5.103}
\| P_{\xi(G_{\alpha}^{N_{i}}, \frac{N_{i}}{4} \leq \cdot \leq 4 N_{i}} u(a_{\alpha, N_{i}}) \|_{L_{x}^{2}(\mathbf{R}^{2})}^{2} \leq \frac{2}{\delta N_{i}} \sum_{J_{l} \subset G_{\alpha}^{N_{i}}} N(J_{l}) \| P_{\xi(a_{l}), \frac{N_{i}}{8} \leq \cdot \leq 8 N_{i}} u(a_{l}) \|_{L_{x}^{2}(\mathbf{R})}^{2}.
\end{equation}

\noindent Therefore,

$$\sum_{1 \leq N_{i} \leq N_{j}} (\frac{N_{i}}{N_{j}}) \sum_{G_{\alpha}^{N_{i}} \cap G_{k}^{N_{j}} \neq \emptyset} \| P_{\xi(G_{\alpha}^{N_{i}}), \frac{N_{i}}{4} \leq \cdot \leq 4 N_{i}} u(a_{\alpha, N_{i}}) \|_{L_{x}^{2}(\mathbf{R})}^{2}$$

$$\leq \frac{1}{N_{j}} \sum_{J_{l} \subset G_{k}^{N_{j}}} \sum_{1 \leq N_{i} \leq N_{j}, \frac{2 N(J_{l})}{\delta} \leq N_{i}} (1 + \frac{N(J_{l})}{\delta}) \| P_{\xi(a_{l}), \frac{N_{i}}{8} \leq \cdot \leq 8 N_{i}} u(a_{l}) \|_{L_{x}^{2}(\mathbf{R})}^{2}$$ $$\lesssim \frac{1}{N_{j}} \sum_{J_{l} \subset G_{k}^{N_{j}}} (1 + \frac{N(J_{l})}{\delta}) \epsilon^{2} \lesssim \epsilon^{2}.$$

\noindent Also,

\begin{equation}\label{5.104}
\sup_{1 \leq N_{i} \leq N_{j}} \sup_{Y_{\alpha'}^{N_{i}} \cap G_{k}^{N_{j}} \neq \emptyset} \| P_{\xi(Y_{\alpha'}^{N_{i}}), \frac{N_{i}}{4} \leq \cdot \leq 4 N_{i}} u(a_{\alpha', N_{i}}) \|_{L_{x}^{2}(\mathbf{R})}^{2} \lesssim \epsilon^{2}.
\end{equation}

\section{Duhamel Terms}
\noindent Now we turn to the Duhamel term

\begin{equation}\label{10.1}
 \int_{a_{\alpha}^{N_{i}}}^{t} e^{i(t - \tau) \Delta} P_{\xi(G_{\alpha}^{N_{i}}), \frac{N_{1}}{4} \leq \cdot \leq 4 N_{1}} (|u(\tau)|^{4} u(\tau)) d\tau,
\end{equation}

\noindent for $N_{1} \geq N_{i}$. We will need the case when $N_{1} >> N_{i}$ for $\S 6$. The arguments for $Y_{\alpha'}^{N_{i}}$ will be virtually identical to the arguments for $G_{\alpha}^{N_{i}}$.

\begin{equation}\label{10.2}
\aligned
 P_{\xi(G_{\alpha}^{N_{i}}), \frac{N_{1}}{4} \leq \cdot \leq 4 N_{1}} (|u(\tau)|^{4} u(\tau)) = \\
O((P_{\xi(G_{\alpha}^{N_{i}}), \geq \frac{N_{1}}{32}} u)(P_{\xi(G_{\alpha}^{N_{i}}), \geq 2^{-10} N_{i}} u) u^{3}) + 5 (P_{\xi(G_{\alpha}^{N_{i}}), \frac{N_{1}}{32} \leq \cdot \leq 32 N_{1}} u)(P_{\xi(G_{\alpha}^{N_{i}}), \leq 2^{-10} N_{i}} u)^{4}.
\endaligned
\end{equation}

\noindent We start with $O((P_{\xi(G_{\alpha}^{N_{i}}), \geq \frac{N_{1}}{32}} u)(P_{\xi(G_{\alpha}^{N_{i}}), \geq 2^{-10} N_{i}} u) u^{3})$.

\begin{theorem}\label{t8.1}
Suppose $N_{1}, N_{2} \geq \frac{N_{i}}{32}$, and $G_{\alpha}^{N_{i}}$ is a green interval.

\begin{equation}\label{8.1}
\aligned
\| |P_{ \xi(G_{\alpha}^{N_{i}}), N_{1}} u| |P_{\xi(G_{\alpha}^{N_{i}}), N_{2}} u| |P_{\xi(G_{\alpha}^{N_{i}}), \leq 2^{-10} N_{i}} u|^{4} \|_{L_{t,x}^{1}(G_{\alpha}^{N_{i}} \times \mathbf{R})} 
\\ \lesssim (\epsilon^{2} + \| u \|_{\tilde{X}_{N_{i}}}^{2}) \| P_{\xi(G_{\alpha}^{N_{i}}), N_{1}} u \|_{U_{\Delta}^{2}(G_{\alpha}^{N_{i}} \times \mathbf{R})}  \| P_{\xi(G_{\alpha}^{N_{i}}), N_{2}} u \|_{U_{\Delta}^{2}(G_{\alpha}^{N_{i}} \times \mathbf{R})} .
\endaligned
\end{equation}
\end{theorem}

\noindent \emph{Proof:} Make a Littlewood - Paley decomposition. Without loss of generality suppose $\xi(G_{\alpha}^{N_{i}}) = 0$.

\begin{equation}\label{8.2}
(\ref{8.1}) \leq \sum_{N_{4} \leq 2^{-10} N_{i}} \| |P_{ N_{1}} u| |P_{N_{2}} u| |P_{\xi(t), N_{4} \leq \cdot \leq 2^{-10} N_{i}} u| |P_{\xi(t), N_{4}} u| |P_{\xi(t), \leq N_{4}} u|^{2} \|_{L_{t,x}^{1}(G_{\alpha}^{N_{i}} \times \mathbf{R})}.
\end{equation}

\noindent Making a bilinear estimate and using $\| u \|_{U_{\Delta}^{2}(J_{l} \times \mathbf{R})} \lesssim_{m_{0}} 1$,

$$\| |P_{N_{1}} u| |P_{N_{2}} u| |P_{\xi(t), N_{4} \leq \cdot \leq 2^{-10} N_{i}} u| |P_{\xi(t), \leq 1} u|^{3}  \|_{L_{t,x}^{1}(J_{l} \times \mathbf{R})}$$

$$\lesssim \| |P_{ N_{1}} u| |P_{\leq 2^{-10} N_{i}} u| \|_{L_{t,x}^{2}(J_{l} \times \mathbf{R})}  \| |P_{ N_{2}} u| |P_{\leq 2^{-10} N_{i}} u| \|_{L_{t,x}^{2}(J_{l} \times \mathbf{R})} \| P_{\xi(t), \leq 1} u \|_{L_{t,x}^{\infty}(J_{l} \times \mathbf{R})}^{2}$$

$$\lesssim \frac{1}{N_{1}^{1/2} N_{2}^{1/2}}  (\| P_{\xi(t), \leq \frac{N(t)}{\delta^{1/2}}} u \|_{L_{t,x}^{\infty}(J_{l} \times \mathbf{R})}^{2} + \| P_{\xi(t), \frac{N(t)}{\delta^{1/2}} \leq \cdot \leq 1} u \|_{L_{t,x}^{\infty}(J_{l} \times \mathbf{R})}^{2})$$  $$\times \| P_{\xi(G_{\alpha}^{N_{i}}), N_{1}} u \|_{U_{\Delta}^{2}(G_{\alpha}^{N_{i}} \times \mathbf{R})}  \| P_{\xi(G_{\alpha}^{N_{i}}), N_{2}} u \|_{U_{\Delta}^{2}(G_{\alpha}^{N_{i}} \times \mathbf{R})} .$$

$$\lesssim \frac{1}{N_{1}^{1/2} N_{2}^{1/2}}  (\frac{N(J_{l})}{\delta^{1/2}} + \frac{N(J_{l})}{\delta} \epsilon^{2})  \| P_{\xi(G_{\alpha}^{N_{i}}), N_{1}} u \|_{U_{\Delta}^{2}(G_{\alpha}^{N_{i}} \times \mathbf{R})}  \| P_{\xi(G_{\alpha}^{N_{i}}), N_{2}} u \|_{U_{\Delta}^{2}(G_{\alpha}^{N_{i}} \times \mathbf{R})} .$$

\noindent Summing over the subintervals $J_{l}$,

\begin{equation}\label{8.3}
\aligned
\| |P_{N_{1}} u| |P_{N_{2}} u| |P_{\xi(t), \leq 2^{-10} N_{i}} u||P_{\xi(t), \leq 1} u|^{3} \|_{L_{t,x}^{1}(G_{\alpha}^{N_{i}} \times \mathbf{R})}  \\ \lesssim (\epsilon^{2} + \| u \|_{\tilde{X}_{N_{i}}}^{2}) \| P_{N_{1}} u \|_{U_{\Delta}^{2}(G_{\alpha}^{N_{i}} \times \mathbf{R})}  \| P_{ N_{2}} u \|_{U_{\Delta}^{2}(G_{\alpha}^{N_{i}} \times \mathbf{R})} .
\endaligned
\end{equation}

\noindent Now we consider the case when $N_{4} \geq 1$. First take the intervals $R_{\beta''}^{N_{4}} \subset G_{\alpha}^{N_{i}}$.

$$\sum_{R_{\beta''}^{N_{4}} \subset G_{\alpha}^{N_{i}}} \sum_{1 \leq N_{4} \leq 2^{-10} N_{i}} \| | P_{N_{1}} u| |P_{N_{2}} u| |P_{\xi(t), N_{4} \leq \cdot \leq 2^{-10} N_{i}} u| |P_{\xi(t), N_{4}} u| |P_{\xi(t), \leq N_{4}} u|^{2} \|_{L_{t,x}^{1}(R_{\beta''}^{N_{4}} \times \mathbf{R})}$$

$$\lesssim \sum_{R_{\beta''}^{N_{4}} \subset G_{\alpha}^{N_{i}}} \sum_{1 \leq N_{4} \leq 2^{-10} N_{i}} \| |P_{N_{1}} u| |P_{\leq 2^{-10} N_{i}} u| \|_{L_{t,x}^{2}(R_{\beta''}^{N_{i}} \times \mathbf{R})}$$ $$\times \| |P_{N_{2}} u| |P_{\leq 2^{-10} N_{i}} u| \|_{L_{t,x}^{2}(R_{\beta''}^{N_{i}} \times \mathbf{R})} \| P_{\xi(t), \leq N_{4}} u \|_{L_{t,x}^{\infty}(R_{\beta''}^{N_{4}} \times \mathbf{R})}^{2}$$

$$\lesssim \frac{1}{N_{1}^{1/2} N_{2}^{1/2}}   \| P_{N_{1}} u \|_{U_{\Delta}^{2}(G_{\alpha}^{N_{i}} \times \mathbf{R})}  \| P_{ N_{2}} u \|_{U_{\Delta}^{2}(G_{\alpha}^{N_{i}} \times \mathbf{R})}  \sum_{1 \leq N_{4} \leq 2^{-10} N_{i}} \sum_{R_{\beta''}^{N_{4}} \subset G_{\alpha}^{N_{i}}} \| P_{\xi(t), \leq N_{4}} u \|_{L_{t,x}^{\infty}(R_{\beta''}^{N_{4}} \times \mathbf{R})}^{2}$$

$$\lesssim \frac{1}{N_{1}^{1/2} N_{2}^{1/2}}   \| P_{N_{1}} u \|_{U_{\Delta}^{2}(G_{\alpha}^{N_{i}} \times \mathbf{R})}  \| P_{ N_{2}} u \|_{U_{\Delta}^{2}(G_{\alpha}^{N_{i}} \times \mathbf{R})}  \sum_{J_{l} \subset G_{\alpha}^{N_{i}}} (\frac{N(J_{l})}{\delta^{1/2}} |\ln(\delta)|	+ \| P_{\xi(t), \frac{N(t)}{\delta^{1/2}} \leq \cdot \leq \frac{N(t)}{\delta}} u \|_{L_{t}^{\infty} L_{x}^{2}(J_{l} \times \mathbf{R})}^{2})$$

$$\lesssim (\epsilon^{2} + \delta^{1/3}) (\frac{N_{i}}{N_{1}^{1/2} N_{2}^{1/2}})  \| P_{N_{1}} u \|_{U_{\Delta}^{2}(G_{\alpha}^{N_{i}} \times \mathbf{R})}  \| P_{ N_{2}} u \|_{U_{\Delta}^{2}(G_{\alpha}^{N_{i}} \times \mathbf{R})}.$$

\noindent Next take the intervals $G_{\beta}^{N_{4}}$. Let $\tilde{G}_{\beta}^{N_{4}} = G_{\beta}^{N_{4}} \cap G_{\alpha}^{N_{i}}$.

$$\sum_{1 \leq N_{4} \leq N_{3} \leq 2^{-10} N_{i}} \sum_{G_{\beta}^{N_{4}} \cap G_{\alpha}^{N_{i}} \neq \emptyset} \| P_{N_{1}} u| |P_{N_{2}} u| |P_{\xi(t), N_{3}} u| |P_{\xi(t), N_{4}} u| |P_{\xi(t), \leq N_{4}} u|^{2} \|_{L_{t,x}^{1}(\tilde{G}_{\beta}^{N_{4}} \times \mathbf{R})}$$

$$\lesssim \sum_{1 \leq N_{4} \leq N_{3} \leq 2^{-10} N_{i}} (\sum_{G_{\beta}^{N_{4}} \cap G_{\alpha}^{N_{i}} \neq \emptyset} \| |P_{N_{1}} u| |P_{\xi(G_{\beta}^{N_{4}}), \frac{N_{4}}{4} \leq \cdot \leq 4 N_{4}} u| \|_{L_{t,x}^{2}(G_{\beta}^{N_{4}} \times \mathbf{R})}^{2})^{1/2}$$

$$\times \| |P_{N_{2}} u| |P_{\xi(t), N_{3}} u | |P_{\xi(t), \leq N_{4}}u|^{2} \|_{L_{t,x}^{2}(G_{\alpha}^{N_{i}} \setminus (\cup R_{\beta''}^{N_{4}}) \times \mathbf{R})}.$$

$$(\sum_{G_{\beta}^{N_{4}} \cap G_{\alpha}^{N_{i}}} \| |P_{N_{1}} u| |P_{\xi(G_{\beta}^{N_{4}}), \frac{N_{4}}{4} \leq \cdot \leq 4 N_{4}} u| \|_{L_{t,x}^{2}(G_{\beta}^{N_{4}} \times \mathbf{R})}^{2})^{1/2}$$ $$\lesssim \frac{1}{N_{1}^{1/2}} \| P_{N_{1}} u \|_{U_{\Delta}^{2}(G_{\alpha}^{N_{i}} \times \mathbf{R})} (\sum_{G_{\beta}^{N_{4}} \cap G_{\alpha}^{N_{i}} \neq \emptyset} \| P_{\xi(G_{\beta}^{N_{4}}), \frac{N_{4}}{4} \leq \cdot \leq 4 N_{4}} u \|_{U_{\Delta}^{2}(G_{\beta}^{N_{4}} \times \mathbf{R})}^{2})^{1/2}.$$

\noindent Making bilinear estimates,

$$ \| |P_{N_{2}} u| |P_{\xi(t), N_{3}} u | |P_{\xi(t), \leq N_{4}} \|_{L_{t,x}^{2}(G_{\alpha}^{N_{i}} \setminus (\cup R_{\beta''}^{N_{4}}) \times \mathbf{R})}$$ $$ \lesssim N_{4} (\sum_{G_{\gamma}^{N_{3}} \cap G_{\alpha}^{N_{i}} \neq \emptyset} \|  |P_{N_{2}} u| |P_{\xi(G_{\gamma}^{N_{3}}), \frac{N_{3}}{4} \leq \cdot \leq 4 N_{3}} u| \|_{L_{t,x}^{2}(G_{\gamma}^{N_{3}} \times \mathbf{R})}^{2})^{1/2}$$ $$+ N_{4}  (\sum_{Y_{\gamma'}^{N_{3}} \cap G_{\alpha}^{N_{i}} \neq \emptyset} \|  |P_{N_{2}} u| |P_{\xi(Y_{\gamma'}^{N_{3}}), \frac{N_{3}}{4} \leq \cdot \leq 4 N_{3}} u| \|_{L_{t,x}^{2}(Y_{\gamma'}^{N_{3}} \times \mathbf{R})}^{2})^{1/2}, $$

\noindent by Sobolev embedding this quantity is

$$\lesssim \frac{N_{4}}{N_{2}^{1/2}} \| P_{N_{2}} u \|_{U_{\Delta}^{2}(G_{\alpha}^{N_{i}} \times \mathbf{R})} (\sum_{G_{\gamma}^{N_{3}}}  \| P_{\xi(G_{\gamma}^{N_{3}}), \frac{N_{3}}{4} \leq \cdot \leq 4 N_{3}} u \|_{U_{\Delta}^{2}(G_{\gamma}^{N_{3}} \times \mathbf{R})}^{2})^{1/2}$$

$$+ \frac{N_{4}}{N_{2}^{1/2}} \| P_{N_{2}} u \|_{U_{\Delta}^{2}(G_{\alpha}^{N_{i}} \times \mathbf{R})} (\sum_{Y_{\gamma'}^{N_{3}}}  \| P_{\xi(Y_{\gamma'}^{N_{3}}), \frac{N_{3}}{4} \leq \cdot \leq 4 N_{3}} u \|_{U_{\Delta}^{2}(Y_{\gamma'}^{N_{3}} \times \mathbf{R})}^{2})^{1/2}.$$

\noindent Again by Cauchy - Schwartz,

$$\sum_{1 \leq N_{4} \leq N_{3} \leq 2^{-10} N_{i}} (\frac{N_{4}}{N_{3}})^{1/2}     ((\frac{N_{4}}{N_{i}})     \sum_{G_{\beta}^{N_{4}} \cap G_{\alpha}^{N_{i}} \neq \emptyset} \| P_{\xi(G_{\beta}^{N_{4}}), \frac{N_{4}}{4} \leq \cdot \leq 4 N_{4}} u \|_{U_{\Delta}^{2}(G_{\beta}^{N_{4}} \times \mathbf{R})}^{2})^{1/2} $$     $$\times  ((\frac{N_{3}}{N_{i}})   \sum_{G_{\gamma}^{N_{3}} \cap G_{\alpha}^{N_{i}} \neq \emptyset} \| P_{\xi(G_{\gamma}^{N_{3}}), \frac{N_{3}}{4} \leq \cdot \leq 4 N_{3}} u \|_{U_{\Delta}^{2}(G_{\beta}^{N_{3}} \times \mathbf{R})}^{2})^{1/2}  \lesssim \| u \|_{\tilde{X}_{N_{i}}}^{2}.$$

\noindent Next,

$$\sum_{1 \leq N_{3} \leq 2^{-10} N_{i}}	 (\frac{N_{3}}{N_{i}}) 		\sharp \{ Y_{\beta'}^{N_{3}} \cap G_{\alpha}^{N_{i}} \neq \emptyset \} 	\leq	\sum_{1 \leq N_{3} \leq 2^{-10} N_{i}} (\frac{N_{3}}{N_{i}}) (\sharp \{ R_{\gamma''}^{N_{3}} \subset G_{\alpha}^{N_{i}} \} + 1)$$

$$\lesssim 1 + \sum_{J_{l} \subset G_{\alpha}^{N_{i}}} \sum_{1 \leq N_{3} \leq \frac{N(J_{l})}{\delta}} \frac{N_{3}}{N_{i}} \lesssim 1.$$

\noindent Because $$\| P_{\xi(Y_{\gamma'}^{N_{3}}), \frac{N_{3}}{4} \leq \cdot \leq 4 N_{3}} u \|_{U_{\Delta}^{2}(Y_{\gamma'}^{N_{3}} \times \mathbf{R})} \lesssim \| u \|_{\tilde{X}_{N_{i}}},$$

$$\sum_{1 \leq N_{3} \leq 2^{-10} N_{i}}	(\frac{N_{3}}{N_{i}})	\sum_{Y_{\gamma'}^{N_{3}} \cap G_{\alpha}^{N_{i}} \neq \emptyset} \| P_{\xi(Y_{\gamma'}^{N_{3}}), \frac{N_{3}}{4} \leq \cdot \leq 4 N_{3}} u \|_{U_{\Delta}^{2}(Y_{\gamma'}^{N_{3}} \times \mathbf{R})}^{2}	\lesssim \| u \|_{\tilde{X}_{N_{i}}}^{2}.$$

\noindent Therefore, by Cauchy-Schwartz,

\begin{equation}\label{8.4}
\aligned
\sum_{1 \leq N_{4} \leq N_{3} \leq 2^{-10} N_{i}} \sum_{G_{\beta}^{N_{4}} \cap G_{\alpha}^{N_{i}} \neq \emptyset} 	\| |P_{N_{1}} u| |P_{N_{2}} u| |P_{\xi(t), N_{3}} u| |P_{\xi(t), N_{4}} u| |P_{\xi(t), \leq N_{4}} u|^{2} \|_{L_{t,x}^{1}(G_{\beta}^{N_{4}} \times \mathbf{R})} \\ \lesssim		(\frac{N_{i}}{N_{1}^{1/2} N_{2}^{1/2}}) \| u \|_{\tilde{X}_{N_{i}}}^{2} \| P_{N_{1}} u \|_{U_{\Delta}^{2}(G_{\alpha}^{N_{i}} \times \mathbf{R})}  \| P_{N_{2}} u \|_{U_{\Delta}^{2}(G_{\alpha}^{N_{i}} \times \mathbf{R})}.
\endaligned
\end{equation}

\noindent Similarly,

\begin{equation}\label{8.5}
\aligned
\sum_{1 \leq N_{4} \leq N_{3} \leq 2^{-10} N_{i}} \sum_{Y_{\beta'}^{N_{4}} \cap G_{\alpha}^{N_{i}} \neq \emptyset} 	\| |P_{N_{1}} u| |P_{N_{2}} u| |P_{\xi(t), N_{3}} u| |P_{\xi(t), N_{4}} u| |P_{\xi(t), \leq N_{4}} u|^{2} \|_{L_{t,x}^{1}(Y_{\beta'}^{N_{4}} \times \mathbf{R})} \\ \lesssim 	(\frac{N_{i}}{N_{1}^{1/2} N_{2}^{1/2}})\| u \|_{\tilde{X}_{N_{i}}}^{2} \| P_{N_{1}} u \|_{U_{\Delta}^{2}(G_{\alpha}^{N_{i}} \times \mathbf{R})} \| P_{N_{2}} u \|_{U_{\Delta}^{2}(G_{\alpha}^{N_{i}} \times \mathbf{R})}.
\endaligned
\end{equation}

\noindent This completes the proof of the theorem. We could make exactly the same arguments for the yellow interval $Y_{\alpha}^{N_{i}}$. $\Box$\vspace{5mm}

\begin{corollary}\label{c8.2}
Making virtually identical arguments,

\begin{equation}\label{8.6}
\aligned
\| |P_{ N_{1}} u| |P_{N_{2}} u| |P_{\leq 2^{-10} N_{i}} (P_{\xi(t), \geq C_{0} N(t)}  u)|^{4} \|_{L_{t,x}^{1}(G_{\alpha}^{N_{i}} \times \mathbf{R})} 	\\
\lesssim (\frac{N_{i}}{N_{1}^{1/2} N_{2}^{1/2}})	(\epsilon^{2} + \| u \|_{\tilde{X}_{N_{i}}}^{2}) \| P_{N_{1}} u \|_{U_{\Delta}^{2}(G_{\alpha}^{N_{i}} \times \mathbf{R})}  \| P_{ N_{2}} u \|_{U_{\Delta}^{2}(G_{\alpha}^{N_{i}} \times \mathbf{R})}	\| u_{\xi(t), \geq C_{0} N(t)} \|_{L_{t}^{\infty} L_{x}^{2}(G_{\alpha}^{N_{i}} \times \mathbf{R})}^{2} .
\endaligned
\end{equation}

\begin{equation}\label{8.7}
\aligned
\| |P_{ N_{1}} u| |P_{N_{2}} u| | P_{\xi(t), \leq C_{0} N(t)}  u|^{4} \|_{L_{t,x}^{1}(G_{\alpha}^{N_{i}} \times \mathbf{R})}  \\
\lesssim C_{0}(\frac{N_{i}}{N_{1}^{1/2} N_{2}^{1/2}}) (\epsilon^{2} + \| u \|_{\tilde{X}_{N_{i}}}^{2}) (\sup_{J_{l} \subset G_{\alpha}^{N_{i}}} \| P_{N_{1}} u \|_{U_{\Delta}^{2}(J_{l} \times \mathbf{R})})	 (\sup_{J_{l} \subset G_{\alpha}^{N_{i}}} \| P_{N_{2}} u \|_{U_{\Delta}^{2}(J_{l} \times \mathbf{R})}).
\endaligned
\end{equation}
\end{corollary}

\begin{theorem}\label{t8.4}
\noindent For $N_{1} \geq N_{i}$,

\begin{equation}\label{8.10}
\aligned
\| 	P_{N_{1}}((P_{\geq \frac{N_{1}}{32}} u)(P_{\geq 2^{-10} N_{i}} u) u^{3}) \|_{DU_{\Delta}^{2}(G_{\alpha}^{N_{i}} \times \mathbf{R})}		\lesssim	\\
(\epsilon^{2} + \| u \|_{\tilde{X}_{N_{i}}}^{2})	\sum_{N_{2} \geq \frac{N_{i}}{32}}	(\frac{N_{i}}{N_{2}})^{1/4}		\| P_{\xi(G_{\alpha}^{N_{i}}), N_{2}} u \|_{U_{\Delta}^{2}(G_{\alpha}^{N_{i}} \times \mathbf{R})}.
\endaligned
\end{equation}
\end{theorem}

\noindent \emph{Proof:}	Take $v$ supported on $|\xi| \sim N_{1}$, $\| v \|_{V_{\Delta}^{2}(G_{\alpha}^{N_{i}} \times \mathbf{R})} 	= 	1$.

$$\int_{G_{\alpha}^{N_{i}}}	\langle v, (P_{> \frac{N_{1}}{32}} u) (P_{> 2^{-10} N_{i}} u) u^{3} \rangle dt	\leq \| |v| |P_{> \frac{N_{1}}{32}} u| |P_{> 2^{-10} N_{i}} u| |u|^{3}	\|_{L_{t,x}^{1}(G_{\alpha}^{N_{i}} \times \mathbf{R})}$$

$$\lesssim	\| |v| |P_{> \frac{N_{1}}{32}} u| |P_{> 2^{-10} N_{i}} u| |P_{\geq 2^{-10} N_{i}} u|^{3} \|_{L_{t,x}^{1}(G_{\alpha}^{N_{i}} \times \mathbf{R})}	+ \| |v| |P_{> \frac{N_{1}}{32}} u| |P_{> 2^{-10} N_{i}} u| |P_{\leq 2^{-10} N_{i}} u|^{3} \|_{L_{t,x}^{1}(G_{\alpha}^{N_{i}} \times \mathbf{R})}.$$

$$ \| |v| |P_{> \frac{N_{1}}{32}} u| |P_{> 2^{-10} N_{i}} u| |P_{\leq 2^{-10} N_{i}} u|^{3} \|_{L_{t,x}^{1}(G_{\alpha}^{N_{i}} \times \mathbf{R})}	$$	$$\lesssim	\| v \|_{L_{t}^{4} L_{x}^{\infty}(G_{\alpha}^{N_{i}} \times \mathbf{R})}	\| P_{> 2^{-10} N_{i}} u \|_{L_{t}^{4} L_{x}^{\infty}(G_{\alpha}^{N_{i}} \times \mathbf{R})}	\| u \|_{L_{t}^{\infty} L_{x}^{2}(G_{\alpha}^{N_{i}} \times \mathbf{R})}	\| |P_{> \frac{N_{1}}{32}} u | |P_{\leq 2^{-10} N_{i}} u|^{2} \|_{L_{t,x}^{2}(G_{\alpha}^{N_{i}} \times \mathbf{R})}.$$

\noindent By theorem $\ref{t8.1}$,

$$\| |P_{> \frac{N_{1}}{32}} u | |P_{\leq 2^{-10} N_{i}} u|^{2} \|_{L_{t,x}^{2}(G_{\alpha}^{N_{i}} \times \mathbf{R})}	\lesssim (\epsilon + \| u \|_{\tilde{X}_{N_{i}}})	\sum_{N_{2} \geq \frac{N_{1}}{32}} (\frac{N_{i}}{N_{2}})^{1/2}	\| P_{\xi(G_{\alpha}^{N_{i}}), N_{2}} u \|_{U_{\Delta}^{2}(G_{\alpha}^{N_{i}} \times \mathbf{R})}.$$

\noindent Therefore,

$$ \| |v| |P_{> \frac{N_{1}}{32}} u| |P_{> 2^{-10} N_{i}} u| |P_{\leq 2^{-10} N_{i}} u|^{3} \|_{L_{t,x}^{1}(G_{\alpha}^{N_{i}} \times \mathbf{R})}	$$

$$\lesssim (\epsilon^{2} + \| u \|_{\tilde{X}_{N_{i}}}^{2})	\sum_{N_{2} \geq \frac{N_{1}}{32}} (\frac{N_{i}}{N_{2}})^{1/2}	\| P_{\xi(G_{\alpha}^{N_{i}}), N_{2}} u \|_{U_{\Delta}^{2}(G_{\alpha}^{N_{i}} \times \mathbf{R})}.$$

\noindent Next, because $V_{\Delta}^{2} \subset U_{\Delta}^{3}$, by $(\ref{2.8})$

$$\| (v) (P_{N_{2}} u) \|_{L_{t,x}^{3}(G_{\alpha}^{N_{i}} \times \mathbf{R})}	\lesssim (\frac{N_{1}}{N_{2}})^{1/4}  \| v \|_{V_{\Delta}^{2}(G_{\alpha}^{N_{i}} \times \mathbf{R})} \| P_{N_{2}} u \|_{U_{\Delta}^{2}(G_{\alpha}^{N_{i}} \times \mathbf{R})}.$$

\noindent Therefore,

$$\| |v| |P_{N_{2}} u|  |P_{> 2^{-10} N_{i}} u|^{4} \|_{L_{t,x}^{1}(G_{\alpha}^{N_{i}} \times \mathbf{R})}$$	$$\lesssim 	 \| P_{> 2^{-10} N_{i}} u \|_{L_{t}^{\infty} L_{x}^{2}(G_{\alpha}^{N_{i}} \times \mathbf{R})}	\sum_{N_{2} \geq \frac{N_{1}}{32}} (\frac{N_{i}}{N_{2}})^{1/4}	\| P_{N_{2}} u \|_{U_{\Delta}^{2}(G_{\alpha}^{N_{i}} \times \mathbf{R})}	\| P_{> 2^{-10} N_{i}} u \|_{L_{t}^{9/2} L_{x}^{18}(G_{\alpha}^{N_{i}} \times \mathbf{R})}^{3}.$$

$$\lesssim (\epsilon^{2}	+ \| u \|_{\tilde{X}_{N_{i}}}^{2})	\sum_{N_{2} \geq \frac{N_{1}}{32}} (\frac{N_{i}}{N_{2}})^{1/4}	\| P_{N_{2}} u \|_{U_{\Delta}^{2}(G_{\alpha}^{N_{i}} \times \mathbf{R})}.$$

\noindent The proof of theorem $\ref{t8.4}$ is complete. $\Box$

\begin{theorem}\label{t8.4.1}
\noindent Suppose $\| u \|_{\tilde{X}_{N_{i}}} \lesssim 1$. Then

\begin{equation}\label{8.10.1}
\aligned
 \| |P_{\geq \frac{N_{1}}{32}} u| |P_{\geq 2^{-10} N_{i}} u| |u|^{3} \|_{DU_{\Delta}^{2}(G_{\alpha}^{N_{i}} \times \mathbf{R})}  \lesssim (\frac{N_{i}}{N_{1}})^{1/2} \| P_{\geq \frac{N_{1}}{32}} u \|_{U_{\Delta}^{2}(G_{\alpha}^{N_{i}} \times \mathbf{R})} \\ \| P_{\geq \frac{N_{1}}{32}} u \|_{U_{\Delta}^{2}(G_{\alpha}^{N_{i}} \times \mathbf{R})}	(\sum_{2^{-10} N_{i} \leq N_{2} \leq 2^{-10} N_{1}} (\frac{N_{2}}{N_{1}})^{1/4} \| P_{\geq N_{2}} u \|_{L_{t}^{\infty} L_{x}^{2}(G_{\alpha}^{N_{i}} \times \mathbf{R})}^{1/2}).
\endaligned
\end{equation}

\end{theorem}

\noindent \emph{Proof:} Take $\| v \|_{V_{\Delta}^{2}(G_{\alpha}^{N_{i}} \times \mathbf{R})} = 1$.

$$\| |v| |P_{\geq \frac{N_{1}}{32}} u| |P_{\geq 2^{-10} N_{i}} u| |P_{\leq 2^{-10} N_{i}} u|^{3} \|_{L_{t,x}^{1}(G_{\alpha}^{N_{i}} \times \mathbf{R})}$$	$$\lesssim	\| |P_{\geq \frac{N_{1}}{32}} u| |P_{\leq 2^{-10} N_{i}} u|^{2} \|_{L_{t,x}^{2}(G_{\alpha}^{N_{i}} \times \mathbf{R})}	\| v \|_{L_{t}^{4} L_{x}^{\infty}(G_{\alpha}^{N_{i}} \times \mathbf{R})} \| P_{\geq 2^{-10} N_{i}} u \|_{L_{t}^{4} L_{x}^{\infty}(G_{\alpha}^{N_{i}} \times \mathbf{R})} \| u \|_{L_{t}^{\infty} L_{x}^{2}(G_{\alpha}^{N_{i}} \times \mathbf{R})},$$

\noindent which by theorem $\ref{t8.1}$, conservation of mass,

$$\lesssim (\frac{N_{i}}{N_{1}})^{1/2} \| P_{\geq \frac{N_{1}}{32}} u \|_{U_{\Delta}^{2}(G_{\alpha}^{N_{i}} \times \mathbf{R})}.$$

\noindent Next,

$$\| |v| |P_{\geq \frac{N_{1}}{32}} u| |P_{\geq 2^{-10} N_{1}} u| |P_{\geq 2^{-10} N_{i}} u|^{3} \|_{L_{t,x}^{1}(G_{\alpha}^{N_{i}} \times \mathbf{R})}$$	$$\lesssim	\| P_{\geq \frac{N_{1}}{32}} u  \|_{L_{t}^{4} L_{x}^{\infty}(G_{\alpha}^{N_{i}} \times \mathbf{R})}	\| v \|_{L_{t}^{4} L_{x}^{\infty}(G_{\alpha}^{N_{i}} \times \mathbf{R})} \| P_{\geq 2^{-10} N_{1}} u \|_{L_{t}^{\infty} L_{x}^{2}(G_{\alpha}^{N_{i}} \times \mathbf{R})} \| P_{\geq 2^{-10} N_{i}} u \|_{L_{t,x}^{6}(G_{\alpha}^{N_{i}} \times \mathbf{R})}^{3}$$

$$\lesssim \| P_{\geq \frac{N_{1}}{32}} u \|_{U_{\Delta}^{2}(G_{\alpha}^{N_{i}} \times \mathbf{R})}	\| P_{\geq 2^{-10} N_{1}} u \|_{L_{t}^{\infty} L_{x}^{2}(G_{\alpha}^{N_{i}} \times \mathbf{R})}.$$

\noindent Finally, for $2^{-10} N_{i} \leq N_{2} \leq 2^{-10} N_{1}$,

$$\| (P_{\geq \frac{N_{1}}{32}} u)(P_{N_{2}} u) \|_{L_{t,x}^{3}(G_{\alpha}^{N_{i}} \times \mathbf{R})}$$

$$\lesssim	\| P_{\geq \frac{N_{1}}{32}} u \|_{L_{t,x}^{6}(G_{\alpha}^{N_{i}} \times \mathbf{R})}^{1/2} \| (P_{\geq \frac{N_{1}}{32}} u)(P_{N_{2}} u) \|_{L_{t,x}^{2}(G_{\alpha}^{N_{i}} \times \mathbf{R})}^{1/2} \| P_{N_{2}} u \|_{L_{t,x}^{\infty}(G_{\alpha}^{N_{i}} \times \mathbf{R})}^{1/2}$$

$$\lesssim (\frac{N_{2}}{N_{1}})^{1/4} \| P_{\geq \frac{N_{1}}{32}} u \|_{U_{\Delta}^{2}(G_{\alpha}^{N_{i}} \times \mathbf{R})}	\| P_{N_{2}} u \|_{L_{t}^{\infty} L_{x}^{2}(G_{\alpha}^{N_{i}} \times \mathbf{R})}^{1/2}.$$

\noindent Therefore,

$$\| |v| |P_{\geq \frac{N_{1}}{32}} u| |P_{2^{-10} N_{i} \leq \cdot \leq 2^{-10} N_{1}} u| |P_{\geq 2^{-10} N_{i}} u|^{3} \|_{L_{t,x}^{1}(G_{\alpha}^{N_{i}} \times \mathbf{R})}$$

$$\lesssim \| v \|_{L_{t,x}^{6}(G_{\alpha}^{N_{i}} \times \mathbf{R})} \| (P_{\geq \frac{N_{1}}{32}} u)(P_{2^{-10} N_{i} \leq \cdot \leq 2^{-10} N_{1}} u) \|_{L_{t,x}^{3}(G_{\alpha}^{N_{i}} \times \mathbf{R})} \| P_{\geq 2^{-10} N_{i}} u \|_{L_{t,x}^{6}(G_{\alpha}^{N_{i}} \times \mathbf{R})}^{3}$$

$$\lesssim	\sum_{2^{-10} N_{i} \leq N_{2} \leq 2^{-10} N_{1}}	(\frac{N_{2}}{N_{1}})^{1/4} \| P_{\geq \frac{N_{1}}{32}} u \|_{U_{\Delta}^{2}(G_{\alpha}^{N_{i}} \times \mathbf{R})}	\| P_{N_{2}} u \|_{L_{t}^{\infty} L_{x}^{2}(G_{\alpha}^{N_{i}} \times \mathbf{R})}^{1/2}.$$

\noindent This completes the proof of the theorem. $\Box$

\begin{theorem}\label{t8.5}
For $N_{1} \geq N_{i}$, $G_{\alpha}^{N_{i}} = [a_{\alpha}^{N_{i}}, b_{\alpha}^{N_{i}}]$, $\xi(G_{\alpha}^{N_{i}}) = 0$,

\begin{equation}\label{8.11}
\| \int_{a_{\alpha}^{N_{i}}}^{t}	e^{i(t - \tau) \Delta}	P_{N_{1}}((P_{N_{2}} u)(P_{\leq 2^{-10} N_{i}} u)^{4}(\tau)) d\tau \|_{U_{\Delta}^{2}(G_{\alpha}^{N_{i}} \times \mathbf{R})}	\lesssim	(\frac{N_{i}}{N_{1}^{1/2} N_{2}^{1/2}}) \| P_{N_{2}} u \|_{U_{\Delta}^{2}(G_{\alpha}^{N_{i}} \times \mathbf{R})}	(\epsilon^{2} + \| u \|_{\tilde{X}_{N_{i}}}^{2}).
\end{equation}
\end{theorem}

\noindent \emph{Proof:} Let $G_{\beta}^{N_{4}} = [a_{\beta}^{N_{4}}, b_{\beta}^{N_{4}}]$, $Y_{\beta'}^{N_{4}} = [a_{\beta'}^{N_{4}}, b_{\beta'}^{N_{4}}]$, $R_{\beta''}^{N_{4}} = [a_{\beta''}^{N_{4}}, b_{\beta''}^{N_{4}}]$. Let

\begin{equation}\label{8.12}
u_{nl}^{G_{\beta}^{N_{4}}, N_{3}}(t) = \int_{a_{\beta}^{N_{4}}}^{t} e^{i(t - \tau) \Delta}		(P_{N_{2}} u)(P_{\xi(\tau), N_{4}} u)(P_{\xi(\tau), N_{3}} u)(P_{\xi(\tau), \leq N_{4}} u)^{2}(\tau) d\tau,
\end{equation}

\begin{equation}\label{8.13}
u_{nl}^{Y_{\beta'}^{N_{4}}, N_{3}}(t) = \int_{a_{\beta}^{N_{4}}}^{t} e^{i(t - \tau) \Delta}	(P_{N_{2}} u)(P_{\xi(\tau), N_{4}} u)(P_{\xi(\tau), N_{3}} u)(P_{\xi(\tau), \leq N_{4}} u)^{2}(\tau) d\tau,	
\end{equation}

\begin{equation}\label{8.14}
u_{nl}^{R_{\beta''}^{N_{4}}}(t) = \int_{a_{\beta}^{N_{4}}}^{t} e^{i(t - \tau) \Delta}	(P_{N_{2}} u)(P_{\xi(\tau), N_{4}} u)(P_{\xi(\tau), N_{4} \leq \cdot \leq 2^{-10} N_{i}} u)(P_{\xi(\tau), \leq N_{4}} u)^{2}(\tau) d\tau.	
\end{equation}

\noindent Then

\begin{equation}\label{8.14.1}
(\ref{8.11}) 	\lesssim	\sum_{1 \leq N_{4} \leq N_{3} \leq 2^{-10} N_{i}}	[\sum_{G_{\beta}^{N_{4}} \cap G_{\alpha}^{N_{i}}}	\| u_{nl}^{G_{\beta}^{N_{4}}, N_{3}}(b_{\beta}^{N_{4}}) \|_{L_{x}^{2}(\mathbf{R})} + \sum_{Y_{\beta'}^{N_{4}} \cap G_{\alpha}^{N_{i}}}	\| u_{nl}^{Y_{\beta'}^{N_{4}}, N_{3}}(b_{\beta'}^{N_{4}}) \|_{L_{x}^{2}(\mathbf{R})}]
\end{equation}

\begin{equation}\label{8.14.2}
+ \sum_{1 \leq N_{4} \leq 2^{-10} N_{i}}		[\sum_{R_{\beta''}^{N_{4}} \subset G_{\alpha}^{N_{i}}} 	\| u_{nl}^{R_{\beta''}^{N_{4}}}(b_{\beta''}^{N_{4}}) \|_{L_{x}^{2}(\mathbf{R})}	+	(\sum_{R_{\beta''}^{N_{4}} \subset G_{\alpha}^{N_{i}}} 	\| u_{nl}^{R_{\beta''}^{N_{4}}}(t) \|_{U_{\Delta}^{2}(R_{\beta''}^{N_{4}})}^{2})^{1/2}]
\end{equation}

\begin{equation}\label{8.14.3}
+ \sum_{1 \leq N_{4} \leq N_{3} \leq 2^{-10} N_{i}}	[(\sum_{G_{\beta}^{N_{4}} \cap G_{\alpha}^{N_{i}}}	\| u_{nl}^{G_{\beta}^{N_{4}}, N_{3}}(b_{\beta}^{N_{4}}) \|_{U_{\Delta}^{2}(G_{\beta}^{N_{4}} \times \mathbf{R})}^{2})^{1/2} + (\sum_{Y_{\beta'}^{N_{4}} \cap G_{\alpha}^{N_{i}}}	\| u_{nl}^{Y_{\beta'}^{N_{4}}, N_{3}}(b_{\beta'}^{N_{4}}) \|_{U_{\Delta}^{2}(Y_{\beta'}^{N_{4}} \times \mathbf{R})}^{2})^{1/2}]
\end{equation}

\begin{equation}\label{8.14.4}
+ \| 	\int_{a_{\alpha}^{N_{i}}}^{t}	e^{i(t - \tau) \Delta}	(P_{N_{2}} u)(P_{\xi(\tau), \leq 1} u)^{3} u(\tau) d\tau \|_{U_{\Delta}^{2}(G_{\alpha}^{N_{i}} \times \mathbf{R})}.
\end{equation}

\noindent Take $\| F \|_{L^{2}(\mathbf{R})}$ supported on $|\xi| \sim N_{1}$. $$\| P_{N_{1}} (\int_{a_{\beta}^{N_{4}}}^{b_{\beta}^{N_{4}}} e^{i(b_{\beta}^{N_{4}} - \tau) \Delta} (P_{N_{2}} u)(P_{\xi(\tau), N_{3}} u)(P_{\xi(\tau), N_{4}} u)(P_{\xi(\tau), \leq N_{4}} u)^{2}(\tau) d\tau) \|_{L_{x}^{2}(\mathbf{R})}$$ $$= \sup_{\| F \|_{L^{2}(\mathbf{R})} = 1} \int_{a_{\beta}^{N_{4}}}^{b_{\beta}^{N_{4}}} \langle F, e^{i(b_{\beta}^{N_{4}} - \tau) \Delta} (P_{N_{2}} u)(P_{\xi(\tau), N_{3}} u)(P_{\xi(\tau), N_{4}} u)(P_{\xi(\tau), \leq N_{4}} u)^{2}(\tau) d\tau \rangle$$

$$ = \int_{a_{\beta}^{N_{4}}}^{b_{\beta}^{N_{4}}}	\langle e^{i(\tau - b_{\beta}^{N_{4}}) \Delta} F, 	(P_{N_{2}} u)(P_{\xi(\tau), N_{3}} u)(P_{\xi(\tau), N_{4}} u)(P_{\xi(\tau), \leq N_{4}} u)^{2} \rangle d\tau$$

$$\lesssim \| (e^{i(\tau - b_{\beta}^{N_{4}}) \Delta} F)(P_{\xi(G_{\beta}^{N_{4}}), \frac{N_{4}}{4} \leq \cdot 4 N_{4}} u) \|_{L_{t,x}^{2}(G_{\beta}^{N_{4}} \times \mathbf{R})}	\| (P_{N_{2}} u)	(P_{\xi(\tau), N_{3}} u)(P_{\xi(\tau), \leq N_{4}} u)^{2}	\|_{L_{t,x}^{2}(G_{\beta}^{N_{4}} \times \mathbf{R})}.$$

\noindent Making a bilinear estimate,

$$\| (e^{i(\tau - b_{\beta}^{N_{4}}) \Delta} F)(P_{\xi(G_{\beta}^{N_{4}}), \frac{N_{4}}{4} \leq \cdot 4 N_{4}} u) \|_{L_{t,x}^{2}(G_{\beta}^{N_{4}} \times \mathbf{R})}	\lesssim \frac{1}{N_{2}^{1/2}} \| P_{\xi(G_{\beta}^{N_{4}}), \frac{N_{4}}{4} \leq \cdot 4 N_{4}} u \|_{U_{\Delta}^{2}(G_{\beta}^{N_{4}} \times \mathbf{R})}.$$

\noindent By Holder's inequality,

$$\sum_{1 \leq N_{4} \leq N_{3} \leq 2^{-10} N_{i}}	\sum_{G_{\beta}^{N_{4}} \cap G_{\alpha}^{N_{i}}}	\| u_{nl}^{G_{\beta}^{N_{4}}, N_{3}}(b_{\beta}^{N_{4}}) \|_{L_{x}^{2}(\mathbf{R})} $$

$$\lesssim	\sum_{1 \leq N_{4} \leq N_{3} \leq 2^{-10} N_{i}}	\frac{1}{N_{2}^{1/2}} (\sum_{G_{\beta}^{N_{4}} \cap G_{\alpha}^{N_{i}} \neq \emptyset}	 \| P_{\xi(G_{\beta}^{N_{4}}), \frac{N_{4}}{4} \leq \cdot 4 N_{4}} u \|_{U_{\Delta}^{2}(G_{\beta}^{N_{4}} \times \mathbf{R})}^{2})^{1/2}$$	$$\times\| (P_{N_{2}} u)	(P_{\xi(\tau), N_{3}} u)(P_{\xi(\tau), \leq N_{4}} u)^{2}	\|_{L_{t,x}^{2}((G_{\alpha}^{N_{i}} \setminus (\cup R_{\beta''}^{N_{4}})) \times \mathbf{R})}.$$

$$\|  (P_{N_{2}} u)	(P_{\xi(\tau), N_{3}} u)(P_{\xi(\tau), \leq N_{4}} u)^{2} \|_{L_{t,x}^{2}(G_{\alpha}^{N_{i}} \setminus (\cup R_{\beta''}^{N_{4}}) \times \mathbf{R})}$$

$$\lesssim (\sum_{G_{\gamma}^{N_{3}} \cap G_{\alpha}^{N_{i}} \neq \emptyset}    \|  (P_{N_{2}} u)	(P_{\xi(\tau), N_{3}} u)(P_{\xi(\tau), \leq N_{4}} u)^{2} \|_{L_{t,x}^{2}(G_{\gamma}^{N_{3}}  \times \mathbf{R})}^{2})^{1/2}$$ $$+ 	(\sum_{Y_{\gamma'}^{N_{3}}}    \|  (P_{N_{2}} u)	(P_{\xi(\tau), N_{3}} u)(P_{\xi(\tau), \leq N_{4}} u)^{2} \|_{L_{t,x}^{2}(G_{\gamma}^{N_{3}}  \times \mathbf{R})}^{2})^{1/2}.$$

$$\| (P_{N_{2}} u)(P_{\xi(\tau), N_{3}} u) \|_{L_{t,x}^{2}(G_{\beta}^{N_{3}} \times \mathbf{R})}	\lesssim	\frac{1}{N_{2}^{1/2}} \| P_{N_{2}} u \|_{U_{\Delta}^{2}(G_{\alpha}^{N_{i}} \times \mathbf{R})} \| P_{\xi(G_{\beta}^{N_{3}}), \frac{N_{3}}{4} \leq \cdot \leq 4 N_{3}} u \|_{U_{\Delta}^{2}(G_{\beta}^{N_{3}} \times \mathbf{R})}.$$

\noindent Therefore, as in the proof of theorem $\ref{t8.1}$,

$$\sum_{1 \leq N_{4} \leq N_{3} \leq 2^{-10} N_{i}}	\sum_{G_{\beta}^{N_{4}} \cap G_{\alpha}^{N_{i}} \neq \emptyset}	\| u_{nl}^{G_{\beta}^{N_{4}}, N_{3}}(b_{\beta}^{N_{4}}) \|_{L_{x}^{2}(\mathbf{R})}$$

$$\lesssim	\| P_{N_{2}} u \|_{U_{\Delta}^{2}(G_{\alpha}^{N_{i}} \times \mathbf{R})} 	(\sum_{1 \leq N_{4} \leq N_{3} \leq 2^{-10} N_{i}}	(\frac{N_{4}}{N_{3}})^{1/2}		((\frac{N_{4}}{N_{2}})	\sum_{G_{\beta}^{N_{4}} \cap G_{\alpha}^{N_{i}} \neq \emptyset}		\| P_{\xi(G_{\beta}^{N_{4}}), \frac{N_{4}}{4} \leq \cdot \leq 4 N_{4}} u \|_{U_{\Delta}^{2}(G_{\beta}^{N_{4}} \times \mathbf{R})}^{2})^{1/2}$$ 	$$\times ((\frac{N_{3}}{N_{2}})	\sum_{G_{\gamma}^{N_{3}} \cap G_{\alpha}^{N_{i}} \neq \emptyset}		\| P_{\xi(G_{\gamma}^{N_{3}}), \frac{N_{3}}{4} \leq \cdot \leq 4 N_{3}} u \|_{U_{\Delta}^{2}(G_{\beta}^{N_{4}} \times \mathbf{R})}^{2})^{1/2})$$

$$+	\| P_{N_{2}} u \|_{U_{\Delta}^{2}(G_{\alpha}^{N_{i}} \times \mathbf{R})} 	(\sum_{1 \leq N_{4} \leq N_{3} \leq 2^{-10} N_{i}}	(\frac{N_{4}}{N_{3}})^{1/2}		((\frac{N_{4}}{N_{2}})	\sum_{G_{\beta}^{N_{4}} \cap G_{\alpha}^{N_{i}} \neq \emptyset}		\| P_{\xi(G_{\beta}^{N_{4}}), \frac{N_{4}}{4} \leq \cdot \leq 4 N_{4}} u \|_{U_{\Delta}^{2}(G_{\beta}^{N_{4}} \times \mathbf{R})}^{2})^{1/2}$$ 	$$\times ((\frac{N_{3}}{N_{2}})	\sum_{Y_{\gamma'}^{N_{3}} \cap G_{\alpha}^{N_{i}} \neq \emptyset}		\| P_{\xi(Y_{\gamma'}^{N_{3}}), \frac{N_{3}}{4} \leq \cdot \leq 4 N_{3}} u \|_{U_{\Delta}^{2}(G_{\gamma}^{N_{3}} \times \mathbf{R})}^{2})^{1/2})$$

$$\lesssim \| P_{N_{2}} u \|_{U_{\Delta}^{2}(G_{\alpha}^{N_{i}} \times \mathbf{R})}	\| u \|_{\tilde{X}_{N_{i}}}^{2}.$$

\noindent Similarly,

$$\sum_{1 \leq N_{4} \leq N_{3} \leq 2^{-10} N_{i}}	\sum_{Y_{\beta'}^{N_{4}} \cap G_{\alpha}^{N_{i}} \neq \emptyset}	\| u_{nl}^{Y_{\beta'}^{N_{4}}, N_{3}}(b_{\beta'}^{N_{4}}) \|_{L_{x}^{2}(\mathbf{R})}	\lesssim \| P_{N_{2}} u \|_{U_{\Delta}^{2}(G_{\alpha}^{N_{i}} \times \mathbf{R})}	\| u \|_{\tilde{X}_{N_{i}}}^{2}.$$

\noindent This takes care of $(\ref{8.14.1})$. Next take $R_{\beta''}^{N_{4}} = [a_{\beta''}^{N_{4}}, b_{\beta''}^{N_{4}}]$.

$$\int_{a_{\beta''}^{N_{4}}}^{b_{\beta''}^{N_{4}}}	\langle e^{i(\tau - b_{\beta''}^{N_{4}}) \Delta} F, 		(P_{N_{2}} u)(P_{\xi(\tau), \leq \cdot 2^{-10} N_{i}} u)(P_{\xi(\tau), N_{4}} u)(P_{\xi(\tau), \leq N_{4}} u)^{2} \rangle d\tau$$

$$\lesssim \| (e^{i(\tau - b_{\beta''}^{N_{4}} \Delta} F)(P_{ \leq 2^{-10} N_{i}} u) \|_{L_{t,x}^{2}(R_{\beta''}^{N_{4}} \times \mathbf{R})}	\| (P_{N_{2}} u)(P_{ \leq 2^{-10} N_{i}} u) \|_{L_{t,x}^{2}(R_{\beta''}^{N_{4}} \times \mathbf{R})}	\| P_{\xi(\tau), \leq N_{4}} u \|_{L_{t,x}^{\infty}(R_{\beta''}^{N_{4}} \times \mathbf{R})}^{2}$$

$$\lesssim \frac{1}{N_{2}}	 \| P_{N_{2}} u \|_{U_{\Delta}^{2}(G_{\alpha}^{N_{i}} \times \mathbf{R})}	\| P_{\xi(\tau), \leq N_{4}} u \|_{L_{t,x}^{\infty}(R_{\beta''}^{N_{4}} \times \mathbf{R})}^{2}.$$

\noindent Therefore,

$$\sum_{1 \leq N_{4} \leq 2^{-10} N_{i}}	\sum_{R_{\beta''}^{N_{4}} \subset G_{\alpha}^{N_{i}}}	\| u_{nl}^{R_{\beta''}^{N_{4}}}(b_{\beta''}^{N_{4}}) \|_{L_{x}^{2}(\mathbf{R})}	\lesssim \frac{1}{N_{2}} 	 \| P_{N_{2}} u \|_{U_{\Delta}^{2}(G_{\alpha}^{N_{i}} \times \mathbf{R})}	\sum_{J_{l} \subset G_{\alpha}^{N_{i}}}	\| P_{\xi(t), \leq \frac{N(J_{l})}{\delta}} u \|_{L_{t,x}^{\infty}(J_{l} \times \mathbf{R})}$$

$$\lesssim \frac{N_{i}}{N_{2}}	 \epsilon^{2}	\| P_{N_{2}} u \|_{U_{\Delta}^{2}(G_{\alpha}^{N_{i}} \times \mathbf{R})}.$$

\noindent Now take $v$ supported on $|\xi| \sim N_{2}$, $\| v \|_{V_{\Delta}^{2}(R_{\beta''}^{N_{2}} \times \mathbf{R})} = 1$.

$$\int_{R_{\beta''}^{N_{4}}}	\langle v, (P_{N_{2}} u)(P_{\xi(t), N_{4} \leq \cdot \leq 2^{-10} N_{i}} u)(P_{\xi(t), N_{4}} u)(P_{\xi(t), \leq N_{4}} u)^{2} \rangle dt	$$	$$\lesssim	\| (P_{N_{2}} u)(P_{\leq 2^{-10} N_{i}} u) \|_{L_{t,x}^{2}(R_{\beta''}^{N_{4}} \times \mathbf{R})}	\| (v) (P_{\leq 2^{-10} N_{i}} u) \|_{L_{t,x}^{5/2}(R_{\beta''}^{N_{4}} \times \mathbf{R})}	\| P_{\xi(t), \leq N_{4}} u \|_{L_{t,x}^{20}(R_{\beta''}^{N_{4}} \times \mathbf{R})}^{2}$$

$$\lesssim	(\frac{N_{4}}{N_{2}})^{7/10} \epsilon^{7/5}	\| P_{N_{2}} u \|_{U_{\Delta}^{2}(G_{\alpha}^{N_{i}} \times \mathbf{R})}	\| P_{\xi(t), \leq N_{4}} u \|_{L_{t,x}^{\infty}(R_{\beta''}^{N_{4}} \times \mathbf{R})}^{7/5}.$$

\noindent We interpolated $\| u \|_{L_{t,x}^{6}(R_{\beta''}^{N_{4}} \times \mathbf{R})} \lesssim 1$ with $\| P_{\xi(t), \leq N_{4}} u \|_{L_{t,x}^{\infty}}$. By Sobolev embedding and conservation of mass,

$$\sum_{1 \leq N_{4} \leq 2^{-10} N_{i}}	\frac{1}{N_{i}^{7/10}}		(\sum_{R_{\beta''}^{N_{4}} \subset G_{\alpha}^{N_{i}}}	\| P_{\xi(t), \leq N_{4}} u \|_{L_{t,x}^{\infty}(R_{\beta''}^{N_{4}} \times \mathbf{R})}^{14/5})^{1/2}$$

$$\lesssim	\sum_{1 \leq N_{4} \leq 2^{-10} N_{i}}	(\frac{N_{4}}{N_{i}})^{1/5}		(\sum_{R_{\beta''}^{N_{4}} \subset G_{\alpha}^{N_{i}}}	\frac{1}{N_{i}} \| P_{\xi(t), \leq N_{4}} u \|_{L_{t,x}^{\infty}(R_{\beta''}^{N_{4}} \times \mathbf{R})})^{1/2}.$$

\noindent By Holders inequality,

$$\lesssim 	(\sum_{1 \leq N_{4} \leq 2^{-10} N_{i}}	\sum_{R_{\beta''}^{N_{4}} \subset G_{\alpha}^{N_{i}}} (\frac{N_{4}}{N_{i}})^{2/5})^{1/2}	(\frac{1}{N_{i}}	\sum_{1 \leq N_{4} \leq 2^{-10} N_{i}}	\| P_{\xi(t), \leq N_{4}} u \|_{L_{t,x}^{\infty}(R_{\beta''}^{N_{4}} \times \mathbf{R})})^{1/2}	\lesssim	(\frac{N_{i}}{N_{2}})^{7/10}	\epsilon.$$

\noindent Therefore,

$$\sum_{1 \leq N_{4} \leq 2^{-10} N_{i}}	(\sum_{R_{\beta''}^{N_{4}} \subset G_{\alpha}^{N_{i}}}	\| u_{nl}^{R_{\beta''}^{N_{4}}}(t) \|_{U_{\Delta}^{2}(R_{\beta''}^{N_{4}} \times \mathbf{R})}^{2})^{1/2}	\lesssim (\frac{N_{i}}{N_{2}})^{7/10}	\epsilon	\| P_{N_{2}} u \|_{U_{\Delta}^{2}(G_{\alpha}^{N_{i}} \times \mathbf{R})}.$$

\noindent This takes care of $(\ref{8.14.2})$.\vspace{5mm}

\noindent Next, for $G_{\beta}^{N_{4}} = [a_{\beta}^{N_{4}}, b_{\beta}^{N_{4}}]$,

$$\int_{a_{\beta}^{N_{4}}}^{b_{\beta}^{N_{4}}}	\langle v, 	(P_{N_{2}} u)(P_{\xi(t), N_{3}} u)(P_{\xi(t), N_{4}} u)(P_{\xi(t), \leq N_{4}} u)^{2} \rangle dt$$

$$\lesssim	\| v (P_{\xi(G_{\beta}^{N_{4}}), \frac{N_{4}}{4} \leq \cdot \leq 4 N_{4}}) \|_{L_{t,x}^{5/2}(G_{\beta}^{N_{4}} \times \mathbf{R})}	\| (P_{N_{2}} u)(P_{\xi(t), N_{3}} u) \|_{L_{t,x}^{2}(G_{\beta}^{N_{4}} \times \mathbf{R})}	\| P_{\xi(t), \leq N_{4}} u \|_{L_{t,x}^{20}(G_{\beta}^{N_{4}} \times \mathbf{R})}^{2}.$$

\noindent By lemma $\ref{l5.0}$, Sobolev embedding,

$$\| P_{\xi(t), \leq N_{4}} u \|_{L_{t,x}^{20}(G_{\beta}^{N_{4}} \times \mathbf{R})}^{2}	\lesssim	N_{4}^{7/10} 	(1 + \| u \|_{\tilde{X}_{N_{i}}})^{2}.$$

\noindent Because $V_{\Delta}^{2} \subset U_{\Delta}^{5/2}$,

$$\| v (P_{\xi(G_{\beta}^{N_{4}}), \frac{N_{4}}{4} \leq \cdot \leq 4 N_{4}}) \|_{L_{t,x}^{5/2}(G_{\beta}^{N_{4}} \times \mathbf{R})}		\lesssim	\frac{1}{N_{2}^{1/5}}	\| P_{\xi(G_{\beta}^{N_{4}}), \frac{N_{4}}{4} \leq \cdot \leq 4 N_{4}} \|_{U_{\Delta}^{2}(G_{\beta}^{N_{4}} \times \mathbf{R})}.$$

\noindent Finally, because $G_{\beta}^{N_{4}}$ overlaps at most two green intervals at level $N_{3}$ and at most two yellow intervals at level $N_{3}$,

$$\| (P_{\xi(t), N_{3}} u)(P_{N_{2}} u) \|_{L_{t,x}^{2}(G_{\beta}^{N_{4}} \times \mathbf{R})}		\lesssim \frac{1}{N_{2}^{1/2}} \| P_{N_{2}} u \|_{U_{\Delta}^{2}(G_{\alpha}^{N_{i}} \times \mathbf{R})}	\| u \|_{\tilde{X}_{N_{i}}}.$$

\noindent Therefore,

$$\sum_{1 \leq N_{4} \leq N_{3} \leq 2^{-10} N_{i}}	(\sum_{G_{\beta}^{N_{4}} \cap G_{\alpha}^{N_{i}}}	\| u_{nl}^{G_{\beta}^{N_{4}}, N_{3}}(t) \|_{U_{\Delta}^{2}(G_{\beta}^{N_{4}} \times \mathbf{R})}^{2})^{1/2} $$

$$\lesssim \| P_{N_{2}} u \|_{U_{\Delta}^{2}(G_{\alpha}^{N_{i}} \times \mathbf{R})} 	\| u \|_{\tilde{X}_{N_{i}}}	\sum_{1 \leq N_{4} \leq N_{3} \leq 2^{-10} N_{i}}	\frac{N_{4}^{1/5}}{N_{2}^{1/5}} 	((\frac{N_{4}}{N_{2}}) \sum_{G_{\beta}^{N_{4}} \cap G_{\alpha}^{N_{i}} \neq \emptyset}	\| 	P_{\xi(G_{\beta}^{N_{4}}), \frac{N_{4}}{4} \leq \cdot \leq 4 N_{4}} u \|_{U_{\Delta}^{2}(G_{\beta}^{N_{4}} \times \mathbf{R})}^{2})^{1/2}$$

$$\lesssim	(\frac{N_{i}}{N_{2}})^{1/2} 		\| P_{N_{2}} u \|_{U_{\Delta}^{2}(G_{\alpha}^{N_{i}} \times \mathbf{R})} 	\| u \|_{\tilde{X}_{N_{i}}}^{2}	\sum_{1 \leq N_{4} \leq N_{3} \leq 2^{-10} N_{i}}	\frac{N_{4}^{1/5}}{N_{2}^{1/5}}		\lesssim 	(\frac{N_{i}}{N_{2}})^{7/10}	\| P_{N_{2}} u \|_{U_{\Delta}^{2}(G_{\alpha}^{N_{i}} \times \mathbf{R})} 	\| u \|_{\tilde{X}_{N_{i}}}^{2}.$$

\noindent Similarly, using $\sharp \{ Y_{\beta'}^{N_{4}} \cap G_{\alpha}^{N_{i}} \} \lesssim \frac{N_{i}}{N_{4}}$,

$$\sum_{1 \leq N_{4} \leq N_{3} \leq 2^{-10} N_{i}}	(\sum_{G_{\beta}^{N_{4}} \cap G_{\alpha}^{N_{i}}}	\| u_{nl}^{Y_{\beta'}^{N_{4}}, N_{3}}(t) \|_{U_{\Delta}^{2}(G_{\beta}^{N_{4}} \times \mathbf{R})}^{2})^{1/2} 		\lesssim (\frac{N_{i}}{N_{2}})^{7/10}	\| P_{N_{2}} u \|_{U_{\Delta}^{2}(G_{\alpha}^{N_{i}} \times \mathbf{R})}	\| u \|_{\tilde{X}_{N_{i}}}^{2}.$$

\noindent This takes care of $(\ref{8.14.3})$. Finally, for $\hat{F}$ supported on $|\xi| \sim N_{2}$, $\| F \|_{L^{2}(\mathbf{R})} = 1$,

$$\int_{a_{l}}^{b_{l}}	\langle e^{i(t - b_{l}) \Delta} F, 	(P_{N_{2}} u)(P_{\leq 2^{-10} N_{i}} u)(P_{\xi(t), \leq 1} u)^{3} \rangle dt$$	$$\lesssim	\| (e^{i(t - b_{l}) \Delta} F)(P_{\leq 2^{-10} N_{i}} u)	\|_{L_{t,x}^{2}(J_{l} \times \mathbf{R})}	\| (P_{N_{2}} u)(P_{\leq 2^{-10} N_{i}} u)	\|_{L_{t,x}^{2}(J_{l} \times \mathbf{R})}$$	$$\times [\| P_{\xi(t), \leq \frac{N(t)}{\delta^{1/2}}} u \|_{L_{t,x}^{\infty}(J_{l} \times \mathbf{R})}^{2}	+ \| P_{\xi(t), \frac{N(t)}{\delta^{1/2}} \leq \cdot \leq 1} u \|_{L_{t,x}^{\infty}(J_{l} \times \mathbf{R})}].$$

$$\lesssim	\frac{N(J_{l})}{N_{2}}	 \epsilon	\| P_{N_{2}} u \|_{U_{\Delta}^{2}(G_{\alpha}^{N_{i}} \times \mathbf{R})}.$$

$$\int_{J_{l}} \langle 		v, (P_{N_{2}} u)(P_{\leq 2^{-10} N_{i}} u)(P_{\xi(t), \leq 1} u)^{3} \rangle dt	$$	$$\lesssim	\| (P_{N_{2}} u)(P_{\leq 2^{-10} N_{i}} u) \|_{L_{t,x}^{2}(J_{l} \times \mathbf{R})}		\| v(P_{\leq 2^{-10} N_{i}} u) \|_{L_{t,x}^{5/2}(J_{l} \times \mathbf{R})}$$	$$\times [\| P_{\xi(t), \leq \frac{N(t)}{\delta^{1/2}}} u \|_{L_{t,x}^{20}(J_{l} \times \mathbf{R})}^{2}	+ 	\| P_{\xi(t), \frac{N(t)}{\delta^{1/2}} \leq \cdot \leq 1} u \|_{L_{t,x}^{20}(J_{l} \times \mathbf{R})}^{2}$$

$$\lesssim	\epsilon^{7/5}	\frac{N(J_{l})}{N_{2}}^{7/10}	\| P_{N_{2}} u \|_{U_{\Delta}^{2}(G_{\alpha}^{N_{i}} \times \mathbf{R})}.$$

\noindent Let	$$u_{nl}^{J_{l}, \leq 1}(t) = \int_{a_{l}}^{t}	e^{i(t - \tau) \Delta} 	(P_{N_{2}} u)(P_{\leq 2^{-10} N_{i}} u)(P_{\xi(\tau), 1} u)(\tau) d\tau.$$

$$\| \int_{a_{\alpha}^{N_{i}}}^{t} e^{i(t - \tau) \Delta} (P_{N_{2}} u)(P_{\leq 2^{-10} N_{i}} u)(P_{\xi(t), \leq 1} u)^{3}(\tau) d\tau \|_{U_{\Delta}^{2}(G_{\alpha}^{N_{i}} \times \mathbf{R})}$$	

$$\lesssim	\sum_{J_{l} \subset G_{\alpha}^{N_{i}}} \| u_{nl}^{J_{l}, \leq 1}(b_{l}) \|_{L_{x}^{2}(\mathbf{R})} + (\sum_{J_{l} \subset G_{\alpha}^{N_{i}}} \| u_{nl}^{J_{l}, \leq 1} \|_{U_{\Delta}^{2}(J_{l} \times \mathbf{R})}^{2})^{1/2}$$	$$\lesssim \epsilon^{7/5}	(\frac{N_{i}}{N_{2}})^{7/10}	\| P_{N_{2}} u \|_{U_{\Delta}^{2}(G_{\alpha}^{N_{i}} \times \mathbf{R})}.$$

\noindent We have finished the proof of theorem $\ref{t8.4}$. $\Box$\vspace{5mm}

\noindent We combine theorems $\ref{t8.4}$ and $\ref{t8.5}$ to estimate the Duhamel terms for $G_{\alpha}^{N_{i}}$. We apply theorem $\ref{t8.4}$ to estimate the first term in $(\ref{10.2})$ and theorem $\ref{t8.5}$ to estimate the second term in $(\ref{10.2})$. The estimates of the Duhamel terms for $Y_{\alpha'}^{N_{i}}$ follow in identical fashion. Therefore, the proof of lemmas $\ref{l5.2}$ and $\ref{l5.3}$, and consequently theorem $\ref{t5.1}$, is complete.

\section{The case when $\int_{0}^{\infty} N(t)^{3} dt < \infty$}

\noindent In this section we prove

\begin{theorem}\label{t6.0}
There does not exist a one dimensional minimal mass blowup solution to $(\ref{0.1})$, $\mu = +1$, with $N(t) \leq 1$, $$\int_{0}^{\infty} N(t)^{3} dt < \infty.$$
\end{theorem}

\noindent To prove this we prove an intermediate theorem.

\begin{theorem}\label{t6.1}
Suppose $u(t,x)$ is a minimal mass blowup solution to $(\ref{0.1})$, $\mu = \pm 1$, with $N(t) \leq 1$ and $$\int_{0}^{\infty} N(t)^{3} dt = \tilde{K} < \infty.$$ Then

\begin{equation}\label{6.0}
\| u(t,x) \|_{L_{t}^{\infty} \dot{H}_{x}^{2}([0, \infty) \times \mathbf{R})} \lesssim_{m_{0}} \tilde{K}^{2}.
\end{equation}
\end{theorem}

\noindent By $(\ref{9.1.18})$ there exists a uniform $K_{0}$ such that if $M$ is any dyadic integer and $[0, T]$ is a compact interval with

\begin{equation}\label{6.1}
\int_{0}^{T} \int |u(t,x)|^{6} dx dt = M \epsilon_{0}^{6},
\end{equation}

\noindent $$\sum_{J_{l} \subset [0, T]} N(J_{l}) = \delta K \leq \delta K_{0}.$$

\noindent After rescaling, $u(t,x) \mapsto \lambda u(\lambda^{2} t, \lambda x)$, $\lambda = \frac{M}{K}$, by theorem $\ref{t5.1}$,

\begin{equation}\label{6.2}
\| u_{\lambda} \|_{\tilde{X}_{M}([0, \frac{T}{\lambda^{2}}] \times \mathbf{R})} \leq C,
\end{equation}

\noindent with C independent of $T$. For $l \geq 5$ let

\begin{equation}\label{6.6}
\mathcal U(2^{l}) = \sup_{T} \| P_{> 2^{l} K_{0}} u \|_{U_{\Delta}^{2}([0, T] \times \mathbf{R})}.
\end{equation}

\noindent By Duhamel's formula

$$\| P_{> 2^{l} K_{0}} u \|_{U_{\Delta}^{2}([0, T] \times \mathbf{R})}	\lesssim_{m_{0}}	\| P_{> 2^{l} K_{0}} u(T) \|_{L_{x}^{2}(\mathbf{R})}$$

$$+ \| P_{> 2^{l} K_{0}} (|u|^{4} u) \|_{DU_{\Delta}^{2}([0, T] \times \mathbf{R})}.$$

\noindent Take $l \geq 5$. By theorem $\ref{t8.5}$,

$$\| \int_{0}^{t} e^{i(t - \tau) \Delta} (P_{> K_{0}} u)(P_{\leq 2^{-10} K_{0}} u)^{4}(\tau) d\tau	\|_{U_{\Delta}^{2}([0, T] \times \mathbf{R})}	\lesssim	\sum_{K_{0} \leq N_{2}}	\frac{K_{0}}{N_{2}}	\lesssim 1.$$

\noindent Splitting the Duhamel term,

$$\| \| |P_{> K_{0}} u| |P_{\geq 2^{-10} K_{0}} u| |u|^{3} \|_{N^{0}([0, T] \times \mathbf{R})}	\lesssim	 \| |P_{> K_{0}} u| |P_{\geq 2^{-10} K_{0}} u|^{4} \|_{L_{t,x}^{6/5} ([0, T] \times \mathbf{R})}$$

$$ + \| |P_{> K_{0}} u| |P_{\geq 2^{-10} K_{0}} u| |P_{\leq 2^{-10} K_{0}} u|^{3}	\|_{L_{t}^{4/3} L_{x}^{1} ([0, T] \times \mathbf{R})}.$$

$$ \| |P_{> K_{0}} u| |P_{\geq 2^{-10} K_{0}} u|^{4} \|_{L_{t,x}^{6/5} ([0, T] \times \mathbf{R})}		\lesssim	\| P_{> K_{0}} u \|_{L_{t, x}^{6} ([0, T] \times \mathbf{R})}	\| P_{\geq 2^{-10} K_{0}} u \|_{L_{t,x}^{6} ([0, T] \times \mathbf{R})}^{4}	\lesssim 1.$$

\noindent We use $$\| P_{\geq 2^{-10} K_{0}} u \|_{L_{t,x}^{6}([0, T] \times \mathbf{R})}	\lesssim \| u_{\lambda} \|_{\tilde{X}_{M}([0, \frac{T}{\lambda^{2}}] \times \mathbf{R})}$$ along with Littlewood-Paley summation and the definition of the $\tilde{X}_{M}$ seminorm. By theorem $\ref{t8.1}$,

$$\| |P_{> K_{0}} u| |P_{\leq 2^{-10} K_{0}} u|^{3} |P_{> 2^{-10} K_{0}} u| \|_{L_{t}^{4/3} L_{x}^{1}([0, T] \times \mathbf{R})}$$	$$\lesssim \| |P_{> K_{0}} u| |P_{\leq 2^{-10} K_{0}} u|^{2} \|_{L_{t,x}^{2}([0, T] \times \mathbf{R})}^{3/2}	\| P_{> 2^{-10} K_{0}} u \|_{L_{t}^{\infty} L_{x}^{2}([0, T] \times \mathbf{R})}^{1/2}	\| P_{\leq 2^{-10} K_{0}} u \|_{L_{t}^{\infty} L_{x}^{2} ([0, T] \times \mathbf{R})}^{1/2}$$

$$+  \| |P_{> K_{0}} u| |P_{\leq 2^{-10} K_{0}} u|^{2} \|_{L_{t,x}^{2}([0, T] \times \mathbf{R})} \| P_{2^{-10} K_{0} \leq \cdot \leq K_{0}} u \|_{L_{t}^{4} L_{x}^{\infty}([0, T] \times \mathbf{R})}	\| u \|_{L_{t}^{\infty} L_{x}^{2}([0, T] \times \mathbf{R})}	\lesssim 1.$$

\noindent Therefore, $\mathcal U(2^{l}) \lesssim 1$ when $l \geq 5$. Because $\sum_{J_{l} \subset [0,T]} N(J_{l}) \leq \delta K_{0}$ for any $T$, $|\xi(t) - \xi(0)| \leq 2^{-20} K_{0}$ for all $t \in [0, \infty)$. $$\sum_{J_{l} \subset [0,T]} N(J_{l}) \leq \delta K_{0}$$ also implies $\lim_{t \rightarrow \infty} N(t) = 0$, which implies $\lim_{t \rightarrow \pm \infty} \| P_{2^{l} K_{0}} u(t) \|_{L_{x}^{2}(\mathbf{R})} = 0$ for $l \geq L_{0}$ for some fixed $L_{0}$. Therefore,

\begin{equation}\label{6.4}
\sup_{T} \| P_{> 2^{l} K_{0}} u \|_{U_{\Delta}^{2}([0, T] \times \mathbf{R})} \lesssim \sup_{T} \| P_{> 2^{l} K_{0}} (|u|^{4} u) \|_{DU_{\Delta}^{2}([0, T] \times \mathbf{R})}.
\end{equation}

\noindent By theorem $\ref{t8.5}$,

\begin{equation}\label{6.5}
 \| P_{> 2^{l} K_{0}} ((P_{> 2^{l - 5} K_{0}} u)(P_{\leq 2^{-10} K_{0}} u)^{4}) \|_{DU_{\Delta}^{2}([0, T] \times \mathbf{R})}	\lesssim 2^{-l/2} \| P_{> 2^{l - 5} K_{0}} u \|_{U_{\Delta}^{2}([0, T] \times \mathbf{R})}.
\end{equation}

\noindent By theorem $\ref{t8.4.1}$,

\begin{equation}\label{6.6}
\aligned
 \| P_{> 2^{l} K_{0}} ((P_{> 2^{l - 5} K_{0}} u)(P_{> 2^{-10} K_{0}} u)u^{3}) \|_{DU_{\Delta}^{2}([0, T] \times \mathbf{R})} \\	\lesssim \| P_{> 2^{l - 5} K_{0}} u \|_{U_{\Delta}^{2}([0, T] \times \mathbf{R})} (\sum_{j = 0}^{l - 5} \frac{2^{j/4}}{2^{l/4}} \| P_{> 2^{j} K_{0}} u \|_{L_{t}^{\infty} L_{x}^{2}([0, \infty) \times \mathbf{R})}^{1/2}).
\endaligned
\end{equation}

\noindent Because

\begin{equation}\label{6.7}
 \sup_{T} \| P_{> 2^{j} K_{0}} u \|_{L_{t}^{\infty} L_{x}^{2}([0, T] \times \mathbf{R})} \rightarrow 0
\end{equation}

\noindent as $j \rightarrow \infty$, there exists $L_{0}$ such that for $l \geq L_{0}$,

\begin{equation}\label{6.8}
 \sup_{T} \| P_{> 2^{l} K_{0}} u \|_{U_{\Delta}^{2}([0, T] \times \mathbf{R})} \leq 2^{-15} \sup_{T} \| P_{> 2^{l - 5} K_{0}} u \|_{U_{\Delta}^{2}([0, T] \times \mathbf{R})}.
\end{equation}

\noindent Therefore, $$\sup_{T} \| P_{> 2^{l} K_{0}} u \|_{U_{\Delta}^{2}([0, T] \times \mathbf{R})} \lesssim_{m_{0}} 2^{-3l}$$

\noindent for $l \geq L_{0}$, which proves $u(t) \in L_{t}^{\infty} \dot{H}_{x}^{2}([0, \infty) \times \mathbf{R})$, and

\begin{equation}\label{6.9}
 \| u(t,x) \|_{L_{t}^{\infty} \dot{H}_{x}^{2}([0, \infty) \times \mathbf{R})} \lesssim_{m_{0}} K_{0}^{2}.
\end{equation}

\noindent $\Box$\vspace{5mm}

\noindent Take some $\eta(t) \rightarrow 0$, possibly very slowly.

\begin{equation}\label{6.18}
\| e^{-ix \cdot \xi(t)} u(t) \|_{\dot{H}^{1}(\mathbf{R})} \lesssim N(t) C(\eta(t)) + \eta(t)^{1/2}.
\end{equation}

\begin{equation}\label{6.19}
E(u(t)) = \frac{1}{2} \int |\nabla u(t,x)|^{2} dx + \frac{1}{6} \int |u(t,x)|^{6} dx.
\end{equation}

\noindent By energy conservation $E(u(t)) = E(u(0))$ for any $t$.\vspace{5mm}

\noindent Now, by Holder's inequality,

$$\int |u(0,x)|^{2} dx \leq \int_{|x - x(0)| \leq \frac{C(\frac{m_{0}^{2}}{1000})}{N(0)}} |u(0,x)|^{2} dx + \frac{m_{0}^{2}}{1000}$$

$$\leq C \| u \|_{L_{x}^{6}(\mathbf{R})}^{2} \frac{C(\frac{m_{0}^{2}}{1000})^{2/3}}{N(0)^{2/3}} + \frac{m_{0}^{2}}{1000}$$

$$\leq C E(u(0))^{1/3} \frac{C(\frac{m_{0}^{2}}{1000})^{2/3}}{N(0)^{2/3}} + \frac{m_{0}^{2}}{1000}.$$

\noindent Now by $(\ref{6.18})$, mass conservation, and the Sobolev embedding theorem, we can choose $t$ sufficiently large so that after a Galilean transformation setting $\xi(t) = 0$,

$$ C E(u(t))^{1/2} \frac{C(\frac{m_{0}^{2}}{1000})}{N(0)} + \frac{m_{0}^{2}}{1000} \leq \frac{m_{0}^{2}}{100}.$$

\noindent But since $E(u(0)) = E(u(t))$, this implies $\int |u(0,x)|^{2} dx \leq \frac{m_{0}^{2}}{100}$, which contradicts mass conservation. This completes the proof of theorem $\ref{t6.0}$. $\Box$\vspace{5mm}

\noindent \textbf{Remark:} We cannot apply these arguments exactly to the focusing case because $E$ is no longer positive definite when $\mu = -1$. These arguments do apply when $\mu = -1$ and $\| u_{0} \|_{L^{2}(\mathbf{R})}$ is less than the mass of the ground state. We will not discuss this matter here.

\section{The case $\int_{0}^{\infty} N(t)^{3} dt = \infty$}
As in the cases when $d \geq 3$, $d = 2$, we defeat this scenario via a frequency localized Morawetz estimate. \cite{CGT1} proved that in the defocusing case

\begin{equation}\label{4.1}
\| u(t,x) \|_{L_{t,x}^{8}([0, T] \times \mathbf{R})}^{8} \lesssim \| u(t) \|_{L_{t}^{\infty} \dot{H}^{1}([0, T] \times \mathbf{R})} \| u(t) \|_{L_{t}^{\infty} L_{x}^{2}([0, T] \times \mathbf{R})}^{3}.
\end{equation}

\noindent See also \cite{PV}. The interaction Morawetz estimate is not positive definite in the focusing case. Let $\chi \in C_{0}^{\infty}(\mathbf{R})$ be an even function,

\begin{equation}\label{4.1.1}
\chi(x) = \left\{
            \begin{array}{ll}
              1, & \hbox{$|x| \leq 1$;} \\
              0, & \hbox{$|x| > 2$.}
            \end{array}
          \right.
\end{equation}

\noindent Here we prove

\begin{theorem}\label{t4.1}
Suppose $u(t,x)$ is a minimal mass blowup solution to $(\ref{0.1})$, $\mu = +1$, on $[0, T]$ with $N(t) \leq 1$,

\begin{equation}\label{4.1.2}
\int_{0}^{T} \int |u(t,x)|^{6} dx dt = M \epsilon_{0}^{6}
\end{equation}

\noindent for some dyadic integer $M$ and for $\| u \|_{L_{t,x}^{6}(J_{l} \times \mathbf{R})} = \epsilon_{0}$,

\begin{equation}\label{4.1.4}
\sum_{J_{l} \subset [0, T]} N(J_{l}) = \delta K.
\end{equation}

\noindent Take $\lambda = \frac{M}{K}$. Let

\begin{equation}\label{4.1.5}
\widehat{Iu}(t,\xi) = \chi(\frac{\xi}{32 M}) \hat{u}_{\lambda}(t,\xi).
\end{equation}

\noindent Then

\begin{equation}\label{4.1.6}
\| Iu_{\lambda} \|_{L_{t,x}^{8}([0, T] \times \mathbf{R})}^{8} \lesssim  o(K)(\frac{M}{K}),
\end{equation}

\noindent $M^{I}(t)$ is a modification of the Morawetz action in \cite{CGT1} (see $(\ref{4.8})$).
\end{theorem}

\noindent \emph{Proof:} Since we are going to work exclusively with the rescaled function $u_{\lambda}$, we will drop the $\lambda$ in our notation and realize that we are working with $u_{\lambda}$ for the rest of this section. \cite{CGT1} defined the action

\begin{equation}\label{4.6}
M(t) = \frac{1}{2} \int_{\mathbf{R}} \int_{\mathbf{R}} a(x - y) |u(t,y)|^{2} Im[\bar{u}(t,x) \partial_{x} u(t,x)] dx dy,
\end{equation}

\begin{equation}\label{4.6.1}
a(x - y) = erf(\frac{x - y}{\epsilon}) = \int_{-\infty}^{\frac{x - y}{\epsilon}} e^{-t^{2}} dt.
\end{equation}

\noindent Taking the limit $\epsilon \rightarrow 0$,

\begin{equation}\label{4.7}
\int_{0}^{T} \int |u(t,x)|^{8} dx dt \lesssim \int_{0}^{T} \partial_{t} M(t) \lesssim \sup_{[0, T]} |M(t)|.
\end{equation}

\noindent Because of conservation of mass and momentum

$$\frac{\partial}{\partial t} \int \int |u(t,y)|^{2} Im[\bar{u}(t,x) \partial_{x} u(t,x)] dx dy = 0,$$

\noindent therefore $$a(x - y) = \int_{0}^{\frac{x - y}{\epsilon}} e^{-t^{2}} dt$$ gives exactly the same Morawetz estimates. We will use this $a(x - y)$ because it is an odd function of $x - y$. Now define the modified action

\begin{equation}\label{4.8}
M_{I}(t) = \frac{1}{2} \int_{\mathbf{R}} \int_{\mathbf{R}} a(x - y) |Iu(t,y)|^{2} Im[\bar{Iu}(t,x) \partial_{x} Iu(t,x)] dx dy.
\end{equation}

\noindent We have

\begin{equation}\label{4.9}
\partial_{t}(Iu) = i \Delta (Iu) - i |Iu|^{4} (Iu) + i |Iu|^{4} (Iu) - i I(|u|^{4} u).
\end{equation}

\noindent If we simply had $$\partial_{t}(Iu) = i \Delta (Iu) - i |Iu|^{4} (Iu),$$ then we would have

\begin{equation}\label{4.10}
\int_{0}^{T} \int |Iu(t,x)|^{8} dx dt \lesssim \int_{0}^{T} \partial_{t} M(t) \lesssim \sup_{[0, T]} |M_{I}(t)|,
\end{equation}

\noindent following the arguments in \cite{CGT} identically. Instead we have

\begin{equation}\label{4.11}
\int_{0}^{T} \int |Iu(t,x)|^{8} dx dt \lesssim \int_{0}^{T} \partial_{t} M_{I}(t) + \mathcal E \lesssim \sup_{[0, T]} |M(t)| + \mathcal E,
\end{equation}

\noindent where

\begin{equation}\label{4.12.1}
\aligned
\mathcal E = \frac{1}{4} \int_{0}^{T} \int \int a(x - y) [I(|u|^{4} \bar{u})(t,y) Iu(t,y) - I(|u|^{4} u)(t,y) \overline{Iu}(t,y)] \\ \times [\overline{Iu}(t,x) \partial_{x} Iu(t,x) - Iu(t,x) \partial_{x} \overline{Iu}(t,x)] dx dy dt
\endaligned
\end{equation}

\begin{equation}\label{4.12.3}
\aligned
+ \frac{1}{4} \int_{0}^{T} \int \int a(x - y) |Iu(t,y)|^{2} [ (|Iu|^{4}(Iu)(t,x) - I(|u|^{4} u)(t,x)) (\partial_{x} \overline{Iu}(t,x))\\ + (|Iu|^{4}(\overline{Iu})(t,x) - I(|u|^{4} \bar{u})(t,x))  (\partial_{x} Iu(t,x)) dx dy dt.
\endaligned
\end{equation}

\begin{equation}\label{4.12.2}
\aligned
+ \frac{1}{4} \int_{0}^{T} \int \int a(x - y) |Iu(t,y)|^{2} [\overline{Iu}(t,x) \partial_{x}(|Iu|^{4}(Iu)(t,x) - I(|u|^{4} u)(t,x)) \\ + Iu(t,x) \partial_{x}(|Iu|^{4}(\overline{Iu})(t,x) - I(|u|^{4} \bar{u})(t,x))] dx dy dt.
\endaligned
\end{equation}

\noindent The interaction Morawetz estimates are Galilean invariant. Indeed, because $a(x - y)$ is an odd function,

\begin{equation}\label{4.13}
\int \int a(x - y) |Iu(t,y)|^{2} Im[i \xi(t) |Iu(t,x)|^{2}] dx dy \equiv 0.
\end{equation}

\noindent Therefore,

\begin{equation}\label{4.13.4}
M_{I}(t) = \int \int a(x - y) |Iu(t,y)|^{2} Im[\overline{Iu}(t,x) (\partial_{x} - i \xi(t)) Iu(t,x)] dx dy.
\end{equation}

\noindent Also, 

$$
 \frac{1}{4} \int_{0}^{T} \int \int a(x - y) [I(|u|^{4} \bar{u})(t,y) Iu(t,y) - I(|u|^{4} u)(t,y) \overline{Iu}(t,y)] (2i \xi(t)) |Iu(t,x)|^{2} dx dy dt$$

$$
+ \frac{1}{4} \int_{0}^{T} \int \int a(x - y) |Iu(t,y)|^{2} [ (|Iu|^{4}(Iu)(t,x) - I(|u|^{4} u)(t,x)) ((-i \xi(t)) \overline{Iu}(t,x))$$	$$+ (|Iu|^{4}(\overline{Iu})(t,x) - I(|u|^{4} \bar{u})(t,x))  ((i \xi(t)) Iu(t,x)) dx dy dt$$

$$
+ \frac{1}{4} \int_{0}^{T} \int \int a(x - y) |Iu(t,y)|^{2} [\overline{Iu}(t,x) (i \xi(t))(|Iu|^{4}(Iu)(t,x) - I(|u|^{4} u)(t,x))$$	$$+ Iu(t,x) (-i \xi(t))(|Iu|^{4}(\overline{Iu})(t,x) - I(|u|^{4} \bar{u})(t,x))] dx dy dt \equiv 0.$$

\noindent Therefore,

\begin{equation}\label{4.13.1}
\aligned
\mathcal E = \frac{1}{4} \int_{0}^{T} \int \int a(x - y) [I(|u|^{4} \bar{u})(t,y) Iu(t,y) - I(|u|^{4} u)(t,y) \overline{Iu}(t,y)] \\ \times [\overline{Iu}(t,x) (\partial_{x} - i \xi(t)) Iu(t,x) - Iu(t,x) (\partial_{x} + i \xi(t)) \overline{Iu}(t,x)] dx dy dt
\endaligned
\end{equation}

\begin{equation}\label{4.13.2}
\aligned
+ \frac{1}{4} \int_{0}^{T} \int \int a(x - y) |Iu(t,y)|^{2} [ (|Iu|^{4}(Iu)(t,x) - I(|u|^{4} u)(t,x)) ((\partial_{x} + i \xi(t)) \overline{Iu}(t,x))\\ + (|Iu|^{4}(\overline{Iu})(t,x) - I(|u|^{4} \bar{u})(t,x))  ((\partial_{x} - i \xi(t)) Iu(t,x)) dx dy dt
\endaligned
\end{equation}

\begin{equation}\label{4.13.3}
\aligned
+ \frac{1}{4} \int_{0}^{T} \int \int a(x - y) |Iu(t,y)|^{2} [\overline{Iu}(t,x) (\partial_{x} - i \xi(t))(|Iu|^{4}(Iu)(t,x) - I(|u|^{4} u)(t,x)) \\ + Iu(t,x) (\partial_{x} + i \xi(t)) (|Iu|^{4}(\overline{Iu})(t,x) - I(|u|^{4} \bar{u})(t,x))] dx dy dt.
\endaligned
\end{equation}

\noindent Let $u_{l} = P_{\leq \frac{M}{32}} u$ and $u_{l} + u_{h} = u$. 

$$|u_{h}|^{2} |u|^{4} \lesssim |u_{h}|^{2} |u_{\leq 2^{-10} M}|^{4} + |u_{h}|^{2} |u_{\geq 2^{-10} M}|^{4}.$$

\noindent By theorem $\ref{t8.1}$, corollary $\ref{c8.2}$,

\begin{equation}\label{4.14}
\| |u_{h}|^{2} |u_{\leq 2^{-10} M}|^{4} \|_{L_{t,x}^{1}([0, T] \times \mathbf{R})}	\lesssim (\sup_{J_{l}} \| P_{\frac{M}{32}} u \|_{U_{\Delta}^{2}(J_{l} \times \mathbf{R})})^{2} C_{0} + \| P_{> N(t) C_{0}} u \|_{L_{t}^{\infty} L_{x}^{2}([0, T] \times \mathbf{R})}^{2}.
\end{equation}

\noindent By Duhamel's formula, $\| u \|_{L_{t}^{4} L_{x}^{\infty}(J_{l} \times \mathbf{R})}	\lesssim_{m_{0}} 1,$ and $N(t) \leq \frac{M}{K}$ on $[0, T]$,

$$\| P_{> \frac{M}{32}} u \|_{U_{\Delta}^{2}(J_{l} \times \mathbf{R})}	\lesssim \| P_{> \frac{M}{32}} u \|_{L_{t}^{\infty} L_{x}^{2}([0, T] \times \mathbf{R})}$$ $$+ \| P_{> 2^{-10} M} u \|_{L_{t}^{\infty} L_{x}^{2}([0, T]  \times \mathbf{R})} \| u \|_{L_{t}^{4} L_{x}^{\infty}([0, T] \times \mathbf{R})} \leq o(1),$$

\noindent with $o(1) \rightarrow 0$ as $K \rightarrow \infty$. Let

\begin{equation}\label{4.15}
C_{0} =(\sup \| P_{> 2^{-10} M} u \|_{U_{\Delta}^{2}(J_{l} \times \mathbf{R})})^{-1},
\end{equation}

\noindent $C_{0} \nearrow \infty$ as $K \rightarrow \infty$, so

$$(\sup_{J_{l}} \| P_{> 2^{-10} M} u \|_{U_{\Delta}^{2}(J_{l} \times \mathbf{R})})^{2} C_{0} + \| P_{\xi(t), \geq C_{0} N(t)} u \|_{L_{t}^{\infty} L_{x}^{2}([-T, T] \times \mathbf{R})}^{2} \leq o(1).$$

$$\| |u_{h}|^{2} |u_{\geq 2^{-10} M}|^{4} \|_{L_{t,x}^{1}([0, T] \times \mathbf{R})}$$ $$	\lesssim \| u_{h} \|_{L_{t}^{5} L_{x}^{10}([0, T] \times \mathbf{R})}^{2} \| u_{\geq 2^{-10} M} \|_{L_{t}^{5} L_{x}^{10}([0, T] \times \mathbf{R})}^{3} \| u_{\geq 2^{-10} M} \|_{L_{t}^{\infty} L_{x}^{2}([0, T] \times \mathbf{R})} \leq o(1).$$

\noindent Now we are ready to estimate $$|M_{I}(t)| + |(\ref{4.13.1})| + |(\ref{4.13.2})| + |(\ref{4.13.3})|.$$

\noindent We start with $|M_{I}(t)|$. Because $N(t) \leq \frac{M}{K}$,

\begin{equation}\label{4.16}
|M_{I}(t)| \lesssim \| u \|_{L_{t}^{\infty} L_{x}^{2}([-T, T] \times \mathbf{R})}^{3} \| (\partial_{x} - i \xi(t)) Iu \|_{L_{t}^{\infty} {L}_{x}^{2}([0, T] \times \mathbf{R})} \lesssim o(K) (\frac{M}{K}).
\end{equation}

\noindent Next we take $(\ref{4.13.2})$. Because $I = 1$ on $|\xi| \leq 32M$, $u_{l}$ is supported on $2^{-5} M$,

$$|Iu_{l}|^{4} (Iu_{l}) - I(|u_{l}|^{4} u_{l}) \equiv 0.$$ Because $(\partial_{x} - i \xi(t)) I \lesssim M$,

$$(\ref{4.13.2}) = \frac{5}{2} \int_{-T}^{T} \int \int a(x - y) |Iu(t,y)|^{2}	Re[[u_{l}^{4} Iu_{h} - I(u_{l}^{4} u_{h})]	(\partial_{x} - i \xi(t)) (Iu)] (t,x) dx dy dt$$

$$ + M \| |u_{h}|^{2} |u|^{4} \|_{L_{t,x}^{1}([-T, T] \times \mathbf{R})} \| Iu \|_{L_{t}^{\infty} L_{x}^{2}([-T, T] \times \mathbf{R})}^{2}.$$

\noindent Also, it suffices for us to consider only $P_{\geq 8M} u$ since we will have cancellation otherwise. Make a Littlewood - Paley decomposition. By the fundamental theorem of calculus,

$$|m(\xi + \xi_{2}) - m(\xi)| \leq |\xi_{2}| \sup |\partial_{x} m(\xi)|.$$

$$\| |u_{l}|^{4} (I P_{> \frac{M}{4}} u) - I(u_{l}^{4} (P_{> \frac{M}{4}} u)) \|_{L_{t}^{6/5} L_{x}^{6/5}([0, T] \times \mathbf{R})}$$  $$ \lesssim \sum_{N_{5} \leq N_{4} \leq N_{3} \leq N_{2} \leq \frac{M}{32}}	(\frac{N_{2}}{M}) \| (u_{h}) (P_{N_{2}} u) \|_{L_{t,x}^{2}([0, T] \times \mathbf{R})}	\| P_{N_{3}} u \|_{L_{t}^{4} L_{x}^{\infty}([0, T] \times \mathbf{R})}$$	$$\times \| P_{N_{4}} u \|_{L_{t}^{12} L_{x}^{3}([0, T] \times \mathbf{R})}		\| P_{N_{5}} u \|_{L_{t,x}^{\infty}([0, T] \times \mathbf{R})}$$

$$\lesssim  \sum_{N_{5} \leq N_{4} \leq N_{3} \leq N_{2} \leq 2^{-10} M}	(\frac{N_{2}}{M})	(\frac{M}{N_{2}})^{1/2}	(\frac{M}{N_{3}})^{1/4}	(\frac{M}{N_{4}})^{1/12}	(\frac{N_{5}}{M})^{1/2}	\lesssim 1.$$

\noindent The second to last inequality follows from lemma $\ref{l5.0}$ and $\| u \|_{\tilde{X}_{M}} \leq C$. Meanwhile,

$$\| (\partial_{x} - i \xi(t)) Iu \|_{L_{t,x}^{6}([0, T] \times \mathbf{R})} \lesssim \sum_{N \leq M} N (\frac{M}{N})^{1/6} o(1) \lesssim o(K) (\frac{M}{K}).$$

\noindent Therefore, $|(\ref{4.13.2})| \lesssim o(K) (\frac{M}{K})$.\vspace{5mm}

\noindent Next, integrating by parts,

$$\int_{0}^{T} \int \int a(x - y) |Iu(t,y)|^{2} [ \overline{Iu}(t,x)] (\partial_{x} - i \xi(t))[|Iu|^{4} (Iu) - I(|u|^{4} u)](t,x) dx dy dt$$

$$ = - \int_{0}^{T} \int \int a(x - y) |Iu(t,y)|^{2} [(\partial_{x} + i \xi(t)) \overline{Iu}(t,x)] [|Iu|^{4} (Iu) - I(|u|^{4} u)](t,x) dx dy dt$$

$$ - \int_{0}^{T} \int \int \partial_{x} a(x - y) \overline{Iu}(t,x) [|Iu|^{4} (Iu) - I(|u|^{4} u)](t,x)	|Iu(t,y)|^{2} dx dy dt.$$

\noindent By Young's inequality, since $\| \partial_{x} a(x - y) \|_{L_{x}^{1}(\mathbf{R})} = 1$,

$$\int_{0}^{T} \int \int \partial_{x} a(x - y) \overline{Iu}(t,x) [|Iu|^{4} (Iu) - I(|u|^{4} u)](t,x) |Iu(t,y)|^{2}	dx dy dt$$		$$\lesssim    	\| I(u_{l}^{4} u_{h}) - u_{l}^{4} (Iu_{h}) \|_{L_{t,x}^{6/5}([0, T] \times \mathbf{R})}	\| Iu \|_{L_{t}^{12} L_{x}^{18}([0, T] \times \mathbf{R})}$$

$$+ \| u_{h} \|_{L_{t}^{4} L_{x}^{\infty}([0, T] \times \mathbf{R})}^{2} \| Iu \|_{L_{t}^{12} L_{x}^{6}([0, T] \times \mathbf{R})}^{6}$$

$$ + \| u_{h} \|_{L_{t}^{5} L_{x}^{10}([0, T] \times \mathbf{R})}^{5} \| Iu \|_{L_{t}^{\infty} L_{x}^{2}([0, T] \times \mathbf{R})}^{3}	\leq o(K) (\frac{M}{K}).$$

\noindent Therefore, $(\ref{4.13.3}) = (\ref{4.13.2}) + o(K) (\frac{M}{K})$, so $(\ref{4.13.3}) \leq o(K) (\frac{M}{K})$.\vspace{5mm}

\noindent Finally we turn to $(\ref{4.13.1})$.

$$I(|u|^{4} u) \overline{Iu} - I(|u|^{4} \bar{u})(Iu)		= [ I(|u|^{4} u) - |Iu|^{4} (Iu) ] \overline{Iu} + [|Iu|^{4} (\overline{Iu}) - I(|u|^{4} \bar{u})] Iu.$$

$$\int_{0}^{T} \int \int |Iu(t,x)| |(\partial_{x} - i \xi(t)) Iu(t,x)| |u_{h}(t,y)|^{2} |u(t,y)|^{4} dx dy dt$$	$$\lesssim M \| Iu(t,x) \|_{L_{t}^{\infty} L_{x}^{2}([0, T] \times \mathbf{R})}^{2} \| |u_{h}(t,y)|^{2} |u(t,y)|^{4} \|_{L_{t,x}^{1}([0, T] \times \mathbf{R})}	\lesssim o(K) (\frac{M}{K}).$$

\noindent Finally, since $$I(u_{l}^{5}) - (Iu_{l})^{5} \equiv 0,$$ it remains to evaluate

$$\int_{0}^{T} \int \int |Iu(t,x)| |(\partial_{x} - i \xi(t)) Iu(t,x)| u_{l}(t,y)^{5} (P_{\geq M} u(t,y)) a(x - y)dx dy dt$$

$$ = \int_{0}^{T} \int \int a(x - y) |Iu(t,x)| |(\partial_{x} - i \xi(t)) Iu(t,x)| \frac{\Delta}{\Delta} [u_{l}(t,y)^{5} (P_{\geq M} u(t,y))] dx dy dt$$

\noindent Integrating by parts

$$ \lesssim \int_{0}^{T} \int \int |Iu(t,x)| |(\partial_{x} - i \xi(t)) Iu(t,x)| (\partial_{x} a(x - y)) \frac{1}{M} u_{l}(t,y)^{5} u_{h}(t,x) dx dy dt.$$

\noindent Again by Young's inequality,

$$\lesssim \| u_{h}(t,y) \|_{L_{t}^{4} L_{x}^{\infty}([0, T] \times \mathbf{R})} \| (\partial_{x} - i \xi(t)) Iu(t,x) \|_{L_{t}^{\infty} L_{x}^{2}([0, T] \times \mathbf{R})}	\| Iu(t,x) \|_{L_{t}^{12} L_{x}^{6}([0, T] \times \mathbf{R})}^{6}	\lesssim o(K) (\frac{M}{K}).$$

\noindent This completes the proof of theorem $\ref{t4.1}$. $\Box$\vspace{5mm}

\noindent \textbf{Remark:} The only properties of $a(x - y)$  that we used in the estimate of $(\ref{4.13.1})$, $(\ref{4.13.2})$, and $(\ref{4.13.3})$ are $a$ is an odd function and there exists a constant $C$ such that

\begin{equation}\label{4.17}
|a(x)| \leq C,
\end{equation}

\noindent and

\begin{equation}\label{4.18}
\| \partial_{x} a(x) \|_{L^{1}(\mathbf{R})} \leq C.
\end{equation}

\noindent Therefore, we have in fact proved

\begin{theorem}\label{t4.3}
 Suppose $a(t,x)$ is an odd function of $x$ for all $t$,

\begin{equation}\label{4.19.1}
 |a(t,x)| \leq C,
\end{equation}

\begin{equation}\label{4.19.2}
\| \partial_{x} a(t,x) \|_{L^{1}(\mathbf{R})} \leq C.
\end{equation}

\noindent Then if $u(t,x)$ is a minimal mass blowup solution to $(\ref{0.1})$, $\mu = \pm 1$,

\begin{equation}\label{4.20.1}
\aligned
 \frac{1}{4} \int_{0}^{T} \int \int a(t, x - y) [I(|u|^{4} \bar{u})(t,y) Iu(t,y) - I(|u|^{4} u)(t,y) \overline{Iu}(t,y)] \\ \times [\overline{Iu}(t,x) (\partial_{x} - i \xi(t)) Iu(t,x) - Iu(t,x) (\partial_{x} + i \xi(t)) \overline{Iu}(t,x)] dx dy dt \lesssim_{m_{0}, d} o(K) C,
\endaligned
\end{equation}

\begin{equation}\label{4.20.2}
\aligned
 \frac{1}{4} \int_{0}^{T} \int \int a(t, x - y) |Iu(t,y)|^{2} [ (|Iu|^{4}(Iu)(t,x) - I(|u|^{4} u)(t,x)) ((\partial_{x} + i \xi(t)) \overline{Iu}(t,x))\\ + (|Iu|^{4}(\overline{Iu})(t,x) - I(|u|^{4} \bar{u})(t,x))  ((\partial_{x} - i \xi(t)) Iu(t,x)) dx dy dt \lesssim_{m_{0}, d} o(K) C,
\endaligned
\end{equation}

\noindent and

\begin{equation}\label{4.20.3}
\aligned
 \frac{1}{4} \int_{0}^{T} \int \int a(x - y) |Iu(t,y)|^{2} [\overline{Iu}(t,x) (\partial_{x} - i \xi(t))(|Iu|^{4}(Iu)(t,x) - I(|u|^{4} u)(t,x)) \\ + Iu(t,x) (\partial_{x} + i \xi(t)) (|Iu|^{4}(\overline{Iu})(t,x) - I(|u|^{4} \bar{u})(t,x))] dx dy dt \lesssim_{m_{0}, d} o(K) C.
\endaligned
\end{equation}

\end{theorem}

\noindent \textbf{Remark:} We will not use the interaction Morawetz estimate of \cite{CGT1}, \cite{PV} for the focusing problem because the interaction Morawetz estimate is not positive definite when $\mu = -1$. Nevertheless, if we did have an interaction Morawetz estimate, theorem $\ref{t4.3}$ implies that the Fourier truncation error is bounded by $o(K)C$ if $a$ satisfies $(\ref{4.19.1})$, $(\ref{4.19.2})$.

\begin{theorem}\label{t4.2}
There does not exist a minimal mass blowup solution with $N(t) \leq 1$, $\int_{0}^{\infty} N(t)^{3} dt = \infty$.
\end{theorem}

\noindent \emph{Proof:} Suppose there did exist a minimal mass blowup solution with $N(t) \leq 1$ and $\int_{0}^{\infty} N(t)^{3} dt = \infty$. Take a compact time interval $[0, T]$ with $$\int_{0}^{T} \int |u(t,x)|^{6} dx dt = M \epsilon_{0}^{6},$$ $M$ a dyadic integer. $[0, T]$ can be partitioned into $M$ small intervals with $\| u(t,x) \|_{L_{t,x}^{6}(J_{l} \times \mathbf{R})} = \epsilon_{0}$. We have $$\sum_{J_{l} \subset [0, T]} N(J_{l}) = \delta K.$$ Rescaling, $u(t,x) \mapsto \lambda^{1/2} u(\lambda^{2} t, \lambda x)$, let $\lambda = \frac{M}{K}$. Let $u_{\lambda}(t,x)$ be the rescaled solution. $[0, \frac{T}{\lambda^{2}}]$ can be partitioned into $M$ small intervals $J_{l}^{\lambda}$, and $$\sum_{J_{l}^{\lambda} \subset [0, \frac{T}{\lambda^{2}}]} N(J_{l}^{\lambda}) = \delta M.$$

\noindent Since $|\xi(t)| \leq 2^{-20} M$ for $t \in [0, T]$, and

\begin{equation}\label{4.21}
\int_{|\xi - \xi(t)| > C(\frac{m_{0}^{2}}{1000}) N(t)} |\hat{u}(t,\xi)|^{2} d\xi,
\end{equation}

\noindent for $K$ sufficiently large,

\begin{equation}\label{4.22}
\frac{m_{0}^{2}}{2} \leq \int_{|x - x(t)| \leq \frac{C(\frac{m_{0}^{2}}{1000})}{N(t)}} |Iu(t,x)|^{2} dx.
\end{equation}

\noindent Therefore,

\begin{equation}\label{4.23}
\int_{0}^{T} N(t)^{3} \frac{m_{0}^{8}}{16} dt \leq \int_{0}^{T} N(t)^{3} (\int_{|x - x(t)| \leq \frac{C(\frac{m_{0}^{2}}{1000})}{N(t)}} |Iu(t,x)|^{2} dx)^{4} dt.
\end{equation}

\noindent By Holder's inequality, and theorem $\ref{t4.1}$,

\begin{equation}\label{4.24}
(\ref{4.23}) \lesssim \int_{0}^{T} N(t)^{3} (\frac{C(\frac{m_{0}^{2}}{1000})}{N(t)})^{3} \| Iu(t) \|_{L_{x}^{8}(\mathbf{R})}^{8} dt \lesssim  o(K) \frac{M}{K}.
\end{equation}

\noindent Since $(\ref{4.23}) \sim M$, the proof of theorem $\ref{t4.2}$ is complete. $\Box$\vspace{5mm}

\noindent This completes the proof of theorem $\ref{t0.2}$.

\newpage

\nocite*
\bibliographystyle{plain}
\bibliography{d=1}

\end{document}